\documentclass[aihp]{imsart}

\RequirePackage{amsthm,amsmath,amsfonts,amssymb}
\RequirePackage[numbers]{natbib}
\RequirePackage[colorlinks,citecolor=blue,urlcolor=blue]{hyperref}

\startlocaldefs
\numberwithin{equation}{section}
\theoremstyle{plain}
\newtheorem{theorem}{Theorem}[section]
\newtheorem{corollary}{Corollary}[section]
\newtheorem{lemma}{Lemma}[section]
\newtheorem{hypothesis}{Hypothesis}[section]
\newtheorem{proposition}{Proposition}[section]
\theoremstyle{remark}
\newtheorem{remark}{Remark}[section]
\newtheorem{definition}{Definition}[section]
\def\R{\mathbb R} 
\def\N{\mathbb N}

\def\E{\mathbb E} 
\def\P{\mathbb P} 
\def\Q{\mathbb Q} 
\def\F{\mathbb F}

\DeclareMathAlphabet{\mathonebb}{U}{bbold}{m}{n} 
\newcommand{\one}{\ensuremath{\mathonebb{1}}} 

\endlocaldefs

\begin{document}

\begin{frontmatter}

\title{Path-dependent 
SDEs with jumps and irregular drift: well-posedness and Dirichlet properties}
\runtitle{Path-dependent SDEs with jumps and irregular drift}

\begin{aug}
\author[A]{\inits{EB}\fnms{Elena}~\snm{Bandini}\ead[label=e1]{elena.bandini7@unibo.it}}
\and
\author[B]{\inits{FR}\fnms{Francesco}~\snm{Russo}\ead[label=e3]{francesco.russo@ensta-paris.fr}}
\address[A]{Dipartimento di Matematica, Piazza di Porta S. Donato 5, 40126 Bologna, Italy\printead[presep={,\ }]{e1}}

\address[B]{Unit\'e de Math\'ematiques Appliqu\'ees, ENSTA Paris, 
828, boulevard des Mar\'echaux, F-91120 Palaiseau, France\printead[presep={,\ }]{e3}}
\end{aug}

\begin{abstract}
We discuss a concept of path-dependent SDE with distributional drift
with possible jumps.
We interpret it via a suitable martingale problem, for which we provide existence and uniqueness.  The corresponding solutions are expected to be Dirichlet processes, nevertheless we give examples of solutions which do not fulfill this property.  In the second part of the paper we indeed state and prove significant new results on the class of Dirichlet processes.
\end{abstract}

\begin{abstract}[language=french]
  Nous introduisons un concept d'EDS d\'ependant de la trajectoire avec drift distributionnel et avec sauts.
  On s'attend que les solutions correspondantes soient  des processus de Dirichlet; n\'eanmoins nous exhibons des
  exemples de solutions ne v\'erifiant pas cette propri\'et\'e. Dans la seconde partie de l'article nous
  prouvons par ailleurs de nouveaux r\'esultats significatifs sur la classe des processus de Dirichlet.
\end{abstract}



\end{frontmatter}

\section{Introduction}

In this paper we discuss path-dependent stochastic differential
equations with possible distributional drift
and jumps of the type
\begin{equation}\label{EFormal}
dX_s =   (\beta'(X_s) +  H(s, X^s))ds+ \sigma(X_s) d W^X_s + k(x) \star (\mu^X- \nu\circ X) + (x- k(x)) \star \mu^X.
\end{equation}
Here $k:\R \rightarrow \R$ is a bounded function  such that   $k(x) =x$ in a neighborhood of $0$,  $\beta:  \R \rightarrow \R$ is a continuous function depending on $k$,  $\sigma:  \R \rightarrow \R$ is a continuous function, not vanishing at zero.
$H:D_{-}(0,T)\rightarrow B(0,T)$ is a bounded and Borel measurable map, where $D_{-}(0,T)$ (resp. $B(0,T)$) will denote  the space of real c\`agl\`ad
(resp. bounded Borel) functions on $[0,\,T]$. 
 $\mu^X(ds\,dx)$ 
is the integer valued random measure on $\R_+ \times \R$ corresponding
 to the jump measure of $X$ and 
  $
  (\nu\circ X)(ds\,dx):= Q(X_{s-},dx)ds,
  $
  where  $Q(\cdot,dx)$ is a   
transition kernel from
  $(\R, \mathcal{B} (\R))$   into $(\R,\mathcal{B} (\R))$, with $Q(y,\{0\})=0$,
such that, for some $\alpha \in [0,\,1]$,
$
y \mapsto \int_{\R } (1 \wedge |x|^{1+\alpha}) \, Q(y,dx)
$
is bounded. 
A solution of \eqref{EFormal} is a couple $(X, \P)$ under which $(\nu\circ X)$ is the compensator of $\mu^X$, $W^X$ is a Brownian motion and $X$ satisfies \eqref{EFormal}. Those solutions will be shown to be not necessarily Dirichlet processes. One of the aim of the paper is indeed to focus on some pathological aspects of Dirichlet processes. 

The Markovian case (with $H=0$) with continuous paths has now a relatively long  history.
Diffusions in the generalized sense were
first considered in the case when the solution is still a semimartingale,
 beginning with \cite{po1}. Later on, many
authors considered special cases of SDEs with generalized coefficients. It is
difficult to quote them all, see for the first contributions \cite{frw1}, \cite{frw2},
\cite{bq}
and \cite{ORT1_PartI} 
for a large bibliography in the semimartingale framework.
 In \cite{frw1} and \cite{frw2}, the authors studied time-independent one-dimensional SDEs of the form
\begin{equation}  \label{SDEone}
\mathrm dX_t = \sigma(X_t) \mathrm dW_t +  \beta'(X_t)\mathrm dt,\quad t\in[0,T],
\end{equation}
whose solutions are possibly non-semimartingale processes, where $\sigma$ is a strictly positive continuous
function and $\beta'$ is the derivative of a real-valued continuous function.
The only supplementary assumption was the
 existence of the function
$\Sigma (x) =  2\int_{0}^{x}\frac{\beta'}{\sigma^{2}}(y)dy$, $x \in \mathbb{R}$,
considered as a suitable limit via regularizations.
Those authors considered solutions in law via the use of a suitable
martingale problem.
The SDE \eqref{SDEone} was also investigated
 by \cite{bq}, where the authors  provided a well-stated framework when $\sigma$
 and $\beta$ are $\gamma$-H\"older
 continuous,  $ \gamma > \frac{1}{2}$. In \cite{rtrut},
 the authors have  also shown that in some cases strong solutions
 exist and pathwise uniqueness holds.
 More recently, in the time-dependent framework (but still one-dimensional),
 a significant contribution
was done by \cite{diel}.
As far as the multidimensional case is concerned,
some important steps were done in \cite{issoglio} and
more recently in \cite{cannizzaro},
when the diffusion matrix is the identity and $\beta'$ is a time-dependent
drift in some proper negative Sobolev space.
In the non-Markovian case, at our knowledge,  the only contribution,
i.e. 
\cite{ORT1_PartI}, refers to the continuous case.
 
 We can find recent significant literature in the Markovian case with L\'evy $\alpha$-stable noise, including the multidimensional case. 
 The first contribution in this direction was one-dimensional and made by \cite{mytnik}.  Further work was done by  \cite{ling}, \cite{deraynal2020}, and \cite{kremp}, the latter even beyond the so-called Young regime. In these works the Brownian motion is replaced by a L\'evy $\alpha$-stable process, which produces
 the regularization by noise.

 Our work includes a non-Markovian drift $H$. Nevertheless,
 even when $H=0$, i.e. in the Markovian case,
 we go in a different direction with respect to the present literature.
 The Markovian component of the generator in our
 case involves local and non-local components. Our equation is driven by a compensated random measure and the regularizing noise is still
 the Brownian motion.
 At our knowledge, our work is the first one in the path-dependent case.
 Our one-dimensional techniques can be adapted to the multidimensional case
 by the use of the Zvonkin transformation, see e.g. \cite{issoglio}.
 We have chosen however to be the most general as possibile in the dimension one:
 in higher dimension the assumptions that one needs are less general.
 
SDEs with distributional drift of the type \eqref{EFormal} will be  interpreted via a suitable martingale problem with respect  to the integro-differential operator  $\mathcal L$ defined in \eqref{Ldistrib}, see Definition \ref{D:mtpb_hom_Mark}. 
This consists  in describing the stochastic behaviour of $f(X)$ under some probability $\P$, 
when $f$ belongs to the  domain $\mathcal D_{\mathcal L}$ defined in \eqref{barD}. In particular, for every $f \in \mathcal D_{\mathcal L}$, $f(X)$ is a special semimartingale. $(X, \P)$ will be a solution of the aforementioned martingale problem. 
$X$ is in general a finite quadratic variation process (i.e. $[X,X]$ exists) but not necessarily a Dirichlet process (i.e. the sum of a martingale and a zero quadratic variation process), see Remark \ref{R:weakDirnew}.
In turn, it  will be shown  to be a weak Dirichlet process. We recall that, given a filtration $\F$,  an $\mathbb F$-weak Dirichlet process is a process of the type 
$X= M+ \Gamma$,
 where $M$ is an $\mathbb F$-local martingale and $\Gamma$ is an $\mathbb F$-orthogonal process vanishing at zero.

Making use of the techniques in \cite{BandiniRusso_RevisedWeakDir}, equation \eqref{EFormal} can be rigorously expressed as 
	\begin{align}\label{decompX}
X&= x_0+ \int_0^\cdot \sigma(X_s) dW^X_s + \int_{]0,\cdot]\times \R} 
k(x)\,(\mu^X(ds\,dx)- Q(X_{s-},dx)d s)
 +\lim_{n \rightarrow \infty}\int_0^\cdot  L f_n(X_{s})ds\notag\\
&+ \int_{]0,\cdot]\times \R} (x - k(x))\mu^X(ds\,dx), 
\end{align}
for every sequence $(f_n) \subseteq \mathcal D_{\mathcal L}$   such that $f_n \underset{n \rightarrow \infty}{\rightarrow} Id$ in $C^{1}$, where
$L$
is the differential operator introduced in \eqref{Lf}
restricted to $\mathcal D_{\mathcal L}.$
The limit appearing in \eqref{decompX} holds in the u.c.p. sense.

We now recall the main results of the paper.
In Section \ref{S21} we provide a suitable definition
for the aforementioned martingale problem, see Definition \ref{D:mtpb_hom_Mark}, 
and state some significant stochastic analysis properties
of a solution.
In particular in Proposition \ref{L: sem_pbmart} we show
that, whenever the drift is a function,
a solution $(X, \P)$ of the classical  martingale problem
is a solution to a Stroock-Varadhan martingale problem with
jumps where the space of test functions is constituted
by $C^2$ bounded functions.
 In Section \ref{S:3.2}, we make use of  a proper  bijective function $h \in \mathcal D_L$ introduced in Proposition \ref{P:equiv}:
 Theorem \ref{equiv_mtgpb_weaksol}
states that  $(X, \P)$ is a solution to the martingale problem if
and only if $(h(X),\P)$ is a semimartingale with given characteristics.
This is a fundamental tool in order to show existence and uniqueness.  
In Proposition \ref{R:Yspecsem} we prove that every solution $X$
is a finite quadratic variation  weak Dirichlet process.
Section \ref{S5}
is devoted to
well-posedness  and continuity properties 
for the martingale problem. Existence and uniqueness is given  in Proposition \ref{P47} in the Markovian case and  in Theorem \ref{T47} in the non-Markovian case. In Proposition \ref{L:Hverified} we study the continuity of the map $\mathcal L$, that is exploited in the companion paper \cite{BandiniRusso_RevisedWeakDir}. Finally, in  Section \ref{S:Dirichlet} we insist  on the fact that the process $X$ is not necessarily a Dirichlet process. 
Moreover, 
we illustrate some new properties related to Dirichlet processes and some pathological aspects.
In Appendix \ref{ATechnical} we justify some technical results, in
Appendix \ref{AFQV}, we discuss the stability of finite quadratic variation processes 
and in Appendix \ref{A:sem} we recall some basic properties of semimartingales with jumps.

\section{Basic notions}

 \subsection{Preliminaries and notations}\label{SPrelim}

 $C^0$ (resp. $C^0_b$) will denote the space of continuous functions (resp. continuous and bounded  functions) on $\R$ equipped with the topology of uniform convergence on each compact (resp. equipped with the topology of uniform convergence).  $C^1$ (resp. $C^2$) will be the space of continuously differentiable (twice continuously differentiable) functions $u:\R\rightarrow \R$.
 They are equipped with the topology of the uniform convergence
 on compact intervals of the functions and the corresponding derivatives.
 $C^{1}_b$  (resp. $C^{2}_b$) is the (topological) intersection of $C^{1}$ and $C^0_b$ (resp. $C^{2}$ and $C^0_b$).
 $D(\R_+)$ will denote the space of real c\`adl\`ag functions on $\R_+$.
We will also indicate by $||\cdot||_{\infty}$ the essential supremum norm and by $||\cdot||_{var}$ the total variation norm.

Let $T>0$ be a finite horizon.
In the following $D(0,\,T)$ (resp.
 $D_{-}(0,\,T)$, $C(0,\,T)$, $B(0,T)$) 
will denote  the space of real c\`adl\`ag
(resp. c\`agl\`ad, continuous, bounded Borel) functions on $[0,\,T]$. Those  spaces    are equipped with the uniform convergence topology.
Given $\eta \in D_{-}(0,\,T)$ 
 we will use the notation 
\begin{align*}
	\eta^{t}(s) := 
	\left\{
	\begin{array}{ll}
	\eta(s) \quad \textup{if}\,\,s <t\\
	\eta(t)\quad \textup{if}\,\,s \geq t.
	\end{array}
	\right.
\end{align*}
For $\eta \in D(0,T)$ we write $\eta^-(t) = \eta(t-)$.

We will denote by $\check \Omega$ the canonical  space, namely
 the space $D(0,\,T)$.
 We will denote  by $\check X$ the canonical process defined by  $\check X_t(\check \omega)= \check \omega(t)$, where $\check \omega$ is a generic element of $\check \Omega$. We also set $\check {\mathcal  F}= \sigma(\check X)$.  Given a topological space $E$, in the sequel $\mathcal{B}(E)$ will denote 
the Borel $\sigma$-field associated with $E$.

A stochastic basis  $(\Omega, \mathcal F, \mathbb F, \P)$ is fixed
throughout the section. We will  suppose that
$\mathbb F$ satisfies the usual conditions.
 Related to   $\mathbb F$,
$\mathcal{P}$ (resp. $\mathcal{\tilde{P}}:=\mathcal{P}\otimes \mathcal{B}(\R)$) will denote the predictable $\sigma$-field on $\Omega \times [0,\,T]$ (resp. on $\tilde{\Omega} := \Omega \times [0,\,T]\times \R$). 

A process $X$ indexed by $\R_+$  will be said to be with integrable variation if the expectation of its total variation is finite. 
$\mathcal{A}$ (resp. $\mathcal{A}_{\textup{loc}}$) will denote  the collection of all adapted processes with   integrable variation (resp.  with locally integrable variation), and    $\mathcal{A}^+$ (resp $\mathcal{A}_{\textup{loc}}^+$)  the collection of all adapted integrable increasing (resp. adapted locally integrable)  processes. 
The significance of locally is the usual one which refers 
to  localization by stopping times, see e.g. (0.39) of  
\cite{jacod_book}.

The concept of random measure   
will be extensively used 	throughout   the paper. 
For a  detailed  discussion on this topic  and the unexplained  notations see 
Chapter I and Chapter II, Section 1, in \cite{JacodBook}, Chapter III in \cite{jacod_book},  and  Chapter XI, Section 1, in \cite{chineseBook}.
In particular, if $\mu$ is a random measure on $[0,\,T]\times \R$, for any measurable real function $H$ defined on $\Omega \times [0,\,T]$, one denotes $H \star \mu_t:= \int_{]0,\,t] \times \R} H(\cdot, s,x) \,\mu(\cdot, ds \,dx)$,
when the stochastic integral in the right-hand side is defined
(with possible infinite values).

We  recall that a transition kernel $Q(a, db)$ of
a measurable space $(A, \mathcal A)$ into another measurable space $(B,\mathcal B)$ is a family
$\{Q(a, \cdot): a \in A\}$ of positive measures on $(B,\mathcal B)$, such that $Q(\cdot ,C)$ is $\mathcal A$-measurable
for each $C \in \mathcal B$, see
for instance in Section 1.1, Chapter I of \cite{JacodBook}.

 Let $X$ be an adapted (c\`adl\`ag) process, so that $X: \Omega \rightarrow \check \Omega$. 
 We set the 
 corresponding  jump measure  $\mu^X$ by
\begin{equation}\label{jumpmeasure}
\mu^X(dt\,dx)= \sum_{s \leq T} \one_{\{\Delta X_s \neq 0\}}\, \delta_{(s, \Delta X_s)}(dt\,dx). 
\end{equation}
We denote by  $\nu^X$ the compensator of $\mu^X$,
         see \cite{JacodBook} (Theorem 1.8, Chapter II).
          From now on for such a process $X$,  $(\mathcal F^X_t)$ will denote  the corresponding canonical filtration, which will be omitted when self-explanatory.

\subsection{Recalls on generators with distributional drift}
Let $\sigma, \beta \in C^0$ such that $\sigma >0$. 
We consider formally the PDE operator of the type 
\begin{equation}\label{Lbeta}
 L
  \psi = \frac{1}{2}\sigma^2  \psi'' + \beta' \psi' 
\end{equation}
in the sense introduced by \cite{frw1, frw2}.
Below  we recall some basic analysis tools coming essentially from Section 2 in \cite{frw1}.

 \begin{definition}
 For a mollifier $\rho$ in the space of Schwartz functions with $\int_\R \rho(x) dx =1$, we set 
 $$
\rho_{\frac{1}{n}}(x) := n \,\phi(n x), \quad  \beta'_n := \beta' \ast \rho_{\frac{1}{n}},  \quad  \sigma_n := \sigma \ast \rho_{\frac{1}{n}}, \quad L_n \psi :=\frac{1}{2} \sigma_n^2 \psi'' + \beta'_n \psi'.
$$
 \end{definition}
\begin{remark}
A priori $\sigma_n, \beta'_n$ and $L_n$ depend on the mollifier $\rho$. 
\end{remark}
In the sequel we will make use of the standing assumption below.
  \begin{hypothesis}\label{H:h}
 We assume the existence of the function 
\begin{equation}\label{Sigma}
\Sigma(x) := \lim_{n \rightarrow \infty} 2 \int_0^x \frac{\beta'_n}{\sigma_n^2} (y) dy
\end{equation}
in $C^0$, independently from the mollifier. 
\end{hypothesis}
  \begin{hypothesis}\label{H:hBIS}
The function $\Sigma$ in \eqref{Sigma}
 is lower bounded,
and 
 $
\int_{-\infty}^0 e^{-\Sigma(x)} dx = \int_{0}^{+\infty} e^{-\Sigma(x)} dx= + \infty.
$
\end{hypothesis}

The following definition and proposition are given in \cite{frw1}, see respectively   Proposition 2.3 and the Definition in Section 2.

 \begin{definition}
 \label{exc}
Set
\begin{align}
	\mathcal D_{L}&:=\{ f \in C^1:\,\,f' e^{\Sigma} \in C^1\}
\label{DLdistr}.
\end{align}
For any $f \in \mathcal D_{L}$, we introduce
 \begin{equation}\label{Lf}
 	 L f= \frac{\sigma^2}{2} (e^{\Sigma} f')' e^{-\Sigma}.
 \end{equation}
This defines without ambiguity $L : \mathcal D_{L} \subset C^1 \rightarrow C^0$, and shows that $f \mapsto Lf$ is a continuous map with respect to the graph topology of $L$, i.e., $L f_n \rightarrow Lf$ in $\mathcal D_L$ if and only if $f_n \rightarrow f$ in $C^1$ and $L f_n \rightarrow Lf$ in $C^0$.
 \end{definition}
  \begin{remark}
 \begin{itemize}
 	\item [(i)]	Setting $\psi=f\in C^1$ in \eqref{Lf}, which does not necessarily belong to $\mathcal D_{L}$ in Definition \ref{exc}, we formally find the expression \eqref{Lbeta}.\item[(ii)] If $f \in \mathcal D_{L}$, \eqref{Lf} is a rigorous representation of \eqref{Lbeta}.
	 \end{itemize}
 \end{remark}
 \begin{proposition}
 \label{P:equiv} Hypothesis \ref{H:h} is equivalent to ask that 
 there is a solution $h \in \mathcal D_L$ to $L h=0$ such that $h(0) =0$ and 
 \begin{equation}\label{h'}
h'(x) := e^{-\Sigma(x)},  \quad x \in \R.
\end{equation}
In particular, $h'(0) =1$, and  $h'$ is strictly positive so that $h$ is bijective and 
	the inverse function  $h^{-1}: \R \rightarrow \R$ is well-defined and continuous. 
 \end{proposition}
 \begin{remark}
   $\mathcal D_{L}$ is a topological subspace of $C^1$, equipped with the graph topology of $L$. Notice that in general the space of smooth functions
   with compact support is not included in  $\mathcal D_{L}$. 
 \end{remark}
\begin{definition} \label{D210}
We denote by  $L^0$   the classical PDE operator 
$
L^0\psi(y) =\frac{\sigma_0^2}{2} \psi''(y)
$
 with 
 		\begin{equation}\label{sigma0}
 		\sigma_0(y)=(\sigma h')(h^{-1}(y)).
 		\end{equation}
	\end{definition}
 
  We recall the following facts, that are collected in Lemma 2.9, and in Propositions 2.10 and  2.13 in \cite{frw1}. 
  \begin{proposition}\label{P:resold}
 Assume Hypotheses \ref{H:h} and \ref{H:hBIS}. The following holds. 
 	\begin{itemize}
 		\item[(a)]
 		 $\mathcal D_{L}$ is dense in  $C^1$.
 \item [(b)]
  For any $f \in \mathcal D_{L}$ we have $f^2 \in \mathcal D_{L}$, and
\begin{equation}\label{Lbetafsquare}
 L f^2= 2 f L f + (f' \sigma)^2.
\end{equation}
In particular, $h^2 \in \mathcal D_{L}$ and $L h^2 = (h'\sigma)^2$. 
 		\item[(c)] 
 	 $\mathcal D_{L^0} = C^2$.
 		 	\item[(d)]
 		 	$\phi \in \mathcal D_{L^0}$ if and only if $\phi \circ h \in \mathcal D_{L}$. Moreover, 
 		$L (\phi \circ h) = (L^0 \,\phi)\circ h$
 		for every $\phi \in C^2$.
 		 	\end{itemize}
 \end{proposition}
We will also need the following assumption referred to some $\alpha \in [0,1]$.
 $C^{1+\alpha}_{\textup{loc}}$ denotes the set of functions belonging to $C^1$ whose derivative belongs to $C^\alpha_{\textup {loc}}$. If $\alpha \in (0,1)$, $C^\alpha_{\textup {loc}}$ denotes  the space of locally $\alpha$-H\"older continuous functions, i.e. the set of functions $f: \R \rightarrow \R$ such that,  for every $M>0$, if $|y| \leq M$, $|z|\leq M$,   there exists $C_M$ such that $|f(y) - f(z)|\leq C_M |y-z|^\alpha$. $C^0_{\textup{loc}}$ (resp. $C^1_{\textup{loc}}$, $C^2_{\textup{loc}}$) denotes by convention  $C^0$ (resp. $C^1$, $C^2$). 
\begin{hypothesis}\label{H:h2}
  The function $\Sigma$ introduced in  \eqref{Sigma} belongs to $C^{\alpha}_{\textup {loc}}$.
\end{hypothesis}

  \begin{remark}\label{Rh}
Hypotheses \ref{H:h}, \ref{H:hBIS} and \ref{H:h2} imply that the function $h$ defined in Proposition \ref{P:equiv} belongs to $C^{1+ \alpha}_{\textup {loc}}$
and  that $h'$  is  bounded.
 \end{remark}

\section{The martingale problem}\label{S21}

\subsection{Formulation of the martingale problem and related properties}\label{S:3.1}
 From here on 
we   fix a   truncation function   $k\in \mathcal K$, where as usual $\mathcal K:=\{k:\R \rightarrow \R \textup{ bounded: } k(x) =x \textup{ in a neighborhood of } 0\}$. 
 Let ${L}$  be a given  operator of the form \eqref{Lbeta} depending on some given functions $\sigma$ and $\beta$.
Assume the validity of Hypotheses \ref{H:h} and \ref{H:hBIS}, 
and let $h$ be the function introduced in Proposition  \ref{P:equiv} related to  $L$.

We will consider  
transition kernels $Q(\cdot,dx)$ from
  $(\R, \mathcal{B} (\R))$   into $(\R,\mathcal{B} (\R))$, with $Q(y,\{0\})=0$,
satisfying the following condition.
	\begin{hypothesis}\label{H:Kmeas}
	For some $\alpha \in [0,\,1]$, 
\begin{align*}
y \mapsto \int_{\R} (1 \wedge |x|^{1+\alpha}) \, Q(y,dx)\quad \textup{is bounded}.
\end{align*}
\end{hypothesis}
\begin{remark}\label{R:Levystable}
  Let $\mu^X(ds\,dx)$ be the jump measure of a L\'evy  $\gamma$-stable
  process with 
  $\gamma=(0,2)$. Then   $\nu^X(ds\,dx) =Q(y,dx)ds$ with $Q(y,dx) =Q_0(dx)= |x|^{-1-\gamma}dx$. In this case, Hypothesis \ref{H:Kmeas} is verified with $\alpha > \gamma-1$. For instance, if $\gamma\in (0,1)$ then $\alpha$ can be chosen to be zero. 
\end{remark}

\begin{remark}\label{R:boundednessTotalVar}
	Hypothesis \ref{H:Kmeas}  means  that, for some $\alpha \in [0,\,1]$, the measure-valued $y \mapsto (1 \wedge |x|^{1+\alpha}) \, Q(y,dx)
$ 
is bounded in the  total variation norm. 
\end{remark}

We consider the topological intersection 
\begin{align}\label{barD}
	\mathcal{D}_{\mathcal L}&:= 
\mathcal{D}_{L} \cap   C^{1+ \alpha}_{\textup{loc}} \cap  C^0_b.
\end{align}
In particular,  $
C^{1+ \alpha}_{\textup{loc}}
\cap  C^0_b$ is a complete metric space 
 equipped with the family of norms 
 \\
 $
 (||f'||_{\alpha, R}
 + ||f||_\infty)_{R\in \N^\ast},
 $
 where 
\begin{equation}\label{normalphaR}
 ||g||_{\alpha, R}
 :=\sup_{x, y: x \neq y, 
|x|\leq R, |y|\leq R}
 \frac{|g(y)-g(x)|}{|y-x|^\alpha}+\sup_{x: |x| \leq R} |g(x)|.
\end{equation}

\begin{proposition}\label{densityDLM}
The set 	$\mathcal{D}_{\mathcal L}$  in \eqref{barD} is dense in $C^1$.
\end{proposition}
\proof
Define the unit partition $\chi: \R \rightarrow \R$ as the smooth function 
\begin{align} \label{chi}
	\chi(a):=\left\{
	\begin{array}{ll}
	1 \quad \textup{if}\,\, a \leq -1\\
	0 \quad \textup{if}\,\, a \geq 0, 
	\end{array}
	\right.
	\end{align} 
	and such that  $\chi(a) \in [0,\,1]$ for $a \in (-1, 0)$. 
	Set 
	\begin{align} \label{chin}
	\chi_N(x):=\chi(|x|- N-1), \quad x \in \R.
		\end{align} 
	Notice that $\chi_N(x)$ 
	is a smooth function and 
	\begin{align*}
	\chi_N(x)
	=\left\{
	\begin{array}{ll}
	1 \quad \textup{if}\,\, |x| \leq N\\
	0 \quad \textup{if}\,\, |x| \geq N+1\\
\in [0,\,1]\quad \textup{otherwise}.
	\end{array}
	\right.
\end{align*}
Let $(\rho_{\frac{1}{N}})$ be a sequence of mollifiers with compact support  converging to the delta measure.
Let $f \in C^1$, and define an approximating  sequence $(f_N)$ of $f$ by setting $f_N(0)=f(0)$ and 
$$
f_N ' := e^{-\Sigma}(f' e^{\Sigma}\chi_N)\ast \rho_{\frac{1}{N}}.
$$
Notice  that $f_N$ is continuous and bounded, being $f_N'$ with  compact support. By Remark \ref{Rh},  since $e^{-\Sigma} \in  C^{\alpha}_{\textup{loc}}$, we get that  $f_N' \in C^{\alpha}_{\textup{loc}}$ and $f_N \in C^{1+\alpha}_{\textup{loc}}$. 
Moreover,  $f_N \in \mathcal D_{L}$ since $f_N' e^\Sigma \in C^{1}$.
Finally, $f_N$ converges to $f$ in $C^1$ since $f_N'$ converges to $f'$ uniformly on compact sets.
\endproof


Consider a functional $H$ defined on $D_{-}(0,T)$
 satisfying the following. 
\begin{hypothesis}\label{H:Hcont}
\begin{enumerate}
	\item $H : D_{-}(0,T) \rightarrow B(0,T)$ is bounded and Borel measurable. 
 \item $H$ fulfills the non-anticipating property, 
 i.e., for every $\eta \in D_{-}(0,\,T)$, 
 $
 H(\eta)(t) = H(\eta^{t})(t)$, $t \in [0,T]$.
\end{enumerate}
\end{hypothesis}
For every $f\in \mathcal{D}_{\mathcal L}$ in \eqref{barD}, we set $\mathcal Lf: D_{-}(0,T)\rightarrow B(0,T)$ as 
%
%
\begin{align}\label{Ldistrib}
( \mathcal L f)(\eta)(t)&:=L f(\eta(t)) + \sigma(\eta(t))H(\eta)(t) f'(\eta(t)) + \int_{\R} (f(\eta(t) + x) -f(\eta(t))
-k(x)\,f'(\eta(t))) Q(\eta(t),dx),
\end{align}
with $L$ the operator defined in \eqref{Lf}. 
%

From here on, for every  $\Phi: D_{-}(0,T)\rightarrow B(0,T)$, 
 we will denote $\Phi(s, \eta):=\Phi(\eta)(s)$, $\eta \in D_{-}(0,T)$, $s \in [0,T]$.
\begin{definition}\label{D:mtpb_hom_Mark}
We say that $(X, \P)$ fulfills the (time-homogeneous) martingale problem with respect to $\mathcal {D}_{\mathcal L}$  in \eqref{barD}, $\mathcal L$ in \eqref{Ldistrib} and $x_0 \in \R$, if for any  $f \in \mathcal {D}_{\mathcal L}$, 
the process
\begin{align}\label{mtg_pb_timehom}
M^f := f(X_{\cdot}) - f(x_0) - \int_0^{\cdot} (\mathcal{L} f)(s,X^{-}) ds
\end{align}
is an $(\mathcal F^X_t)$-local martingale  under $\P$.
\end{definition}
\begin{remark}\label{R:app1}
 Hypothesis \ref{H:Kmeas} implies that 
$
y \mapsto \int_{\R} (1 \wedge |x|^{2}) \, Q(y,dx)
$
is  bounded. 
In particular  
the pair $(X,\P)$ in Definition \ref{D:mtpb_hom_Mark} satisfies  
\begin{equation}\label{int_small_jumps}
		\sum_{s \leq \cdot} |\Delta X_s|^2 
		< \infty \,\,\,\,\textup{a.s.},
\end{equation}
see Proposition C.1 in \cite{BandiniRusso_RevisedWeakDir}.
\end{remark}

\begin{remark}
	Let $k \in \mathcal K$ be a generic truncation function. While $\mathcal D_{\mathcal L}$ in \eqref{barD} does not depend on $k$, clearly $\mathcal L$ defined in \eqref{Ldistrib} a priori depends on $k$, namely $\mathcal L = \mathcal L^k$. In order to formulate a coherent definition, we should allow $\beta$ also depending on $k$, as we will explain below. This in particular forces $L=L^k$ to depend on $k$ as well. 
	
Indeed, let $\tilde k \in \mathcal K$. By \eqref{Ldistrib}, for every $\eta \in D_{-}(0,T)$,  we have 
\begin{align*}
&( \mathcal L^k f)(\eta)(t)-( \mathcal L^{\tilde k} f)(\eta)(t)= L^k f(\eta(t))-L^{\tilde k} f(\eta(t)) + f'(\eta(t)))\int_{\R} (\tilde k(x)
-k(x))\, Q(\eta(t),dx). 
\end{align*}
Let $(X, \P)$ fulfilling the  martingale problem with respect to $\mathcal {D}_{\mathcal L}$  in \eqref{barD}, $\mathcal L^k$ in \eqref{Ldistrib} and $x_0 \in \R$. Then $(X, \P)$ fulfills the  martingale problem with respect to $\mathcal {D}_{\mathcal L}$  in \eqref{barD}, $\mathcal L^{\tilde k}$ in and $x_0$, 
  if and only if 
$$
 \int_0^\cdot (L^k f(X_{s-})-L^{\tilde k}f(X_{s-})) ds = \int_0^\cdot f'(X_{s-})\int_{\R} (k(x)
-\tilde k(x))\, Q(X_{s-},dx)ds.
$$ 
This condition is verified if 
$$
\beta^k(X_{s-}) - \beta^{\tilde k}(X_{s-}) = \int_{\R} (k(x)
-\tilde k(x))\, Q(X_{s-},dx)ds.
$$
With this choice, $L^k$ coincides with $L^{\tilde k}$ and  consequently $\mathcal L^k$ coincides  $\mathcal L^{\tilde k}$. 
\end{remark}

When $\beta'$ is a continuous function,  we  recover the classical martingale problem in the sense of Jacod-Shiryaev, see Proposition \ref{L: sem_pbmart} below. 
In the following $s(\mathcal H, X|\P_{\mathcal H}; B, C, \nu)$ denotes the set of all  solutions $\P$ related to a given  probability $\P_{\mathcal H}$ and  characteristics $(B, C, \nu)$,  see Definition \ref{D:mrtg_pb_Jacod}. 
\begin{proposition}\label{L: sem_pbmart}
Let 
$b, \sigma$ be continuous functions, and set 
\begin{equation}\label{prop_BCnu}
 B_t = \int_0^t (b(\check X_s) + \sigma(\check X_s)\normalcolor H(s,\check X^{-})) ds, \quad C_t = \int_0^t \sigma^2(\check X_s) ds, \quad \nu(ds\,dx)= Q(\check X_s, dx)ds,
\end{equation}
with $Q$ satisfying  Hypothesis \ref{H:Kmeas} with $\alpha=1$, and $H$ satisfying Hypothesis \ref{H:Hcont}. Let $L$  be the operator of the form \eqref{Lbeta}    with 
$\beta' := b$.  Set
$ \mathcal H=\{A \in \mathcal F: \exists A_0 \in \mathcal B(\R) \textup{ such that } A= \{\omega \in \Omega: \omega(0) \in A_0\}\}$ and $\P_{\mathcal H} $ corresponds   to $ \delta_{x_0}$  in the sense that, for any $A \in \mathcal F$, $\P_{\mathcal H}(A) = \delta_{x_0}(A_0)$ with $A_0=\{\omega(0) \in \R: \,\, \omega \in A\}$.


 Then   $\P$ belongs to  $s(\mathcal H, X|\P_{\mathcal H}; B, C, \nu)$  if and only if  $(X, \P)$ is  a  solution to the martingale problem in Definition \ref{D:mtpb_hom_Mark} related to  $ \mathcal { D_\mathcal L}=C^2_b$,   $x_0=\check X_0$ and $\mathcal L$ in 
\eqref{Ldistrib}.
	 \end{proposition}
\proof
By  Theorem \ref{T: equiv_mtgpb_semimart} together with Definition \ref{D:mrtg_pb_Jacod},
$\P$ belongs to  $s(\mathcal H, X|\P_{\mathcal H}; B, C, \nu)$  if and only if, for any $f \in C^2_b$, 
	\begin{align*}
 	f(\check X_t) &- f(\check X_0) - \int_0^t\int_{\R} [f(\check X_{s-} + x) -f(\check X_{s-})-k(x)f'(\check X_{s-})]\,\nu(ds\,dx)\notag\\
 	& - \int_0^t \Big[(b(\check X_s)+\sigma(\check X_s)H(s,\check X^{-})) f'(\check X_{s}) 
 	+\frac{1}{2} \sigma^2(\check X_s) f''(\check X_{s})\Big] ds 
 	\end{align*}
 	is a $\P$-local martingale. This agrees in particular with Definition \ref{D:mtpb_hom_Mark} related to $\mathcal L$ in \eqref{Ldistrib},  $
 \mathcal { D_\mathcal L}= C^2_b$ and to  $x_0=\check X_0$, where $L$ is the operator of the form \eqref{Lbeta}    with 
$\beta' := b$.
\endproof

\subsection{About  equivalent formulations for the martingale problem
}\label{S:3.2}
We provide  an equivalent martingale formulation for $Y= h(X)$, with $h$ the function introduced in  Proposition \ref{P:equiv}. 
 This principle can be extended to general  bijective $C^1$-type transformations.
For any $y \in \R$,  introduce 
\begin{align}
F(y, A)&:=\int_{\R}\one_{A} \, (h(h^{-1}(y)+ w)- h(h^{-1}(y))) \, Q(h^{-1}(y), dw)\label{identity_nu_nuY_F},\quad A \subseteq \R,\\
b(y) &:= (h' \circ h^{-1})(y)\int_{\R}[(h^{-1})'(y)\,k(z)-k(h^{-1}(y + z)- h^{-1}(y))]\,F(y, dz)\label{b(y)}.
\end{align}
For any $\phi \in C^2_b$, we also define 
\begin{align}
 {\bar L}\phi&:= {L^0}\phi +  b \, \phi',\label{tildeL}\\
 \bar H(t, \eta)&:= H(t,h^{-1}(\eta)),\label{barH}\\
(\bar {\mathcal{L}} \phi)(t, \eta)
&:= \bar L \phi(\eta(t)) +\sigma_0(\eta(t))
 \bar H(t,\eta)\normalcolor\phi'(\eta(t))\notag\\
&+  \int_{\R} (\phi(\eta_t + z) -\phi(\eta_t)
-k(z)\,\phi'(\eta(t))\,) F(\eta(t),dz), \quad\eta \in D_{-}(0,\,T),\label{barL}
\end{align} 
with  $L^0$ the  operator in Definition  \ref{D210} and $\sigma_0$ in \eqref{sigma0}.

 \begin{remark}\label{D:2new_mrtg_pb}
 Let  $Q(\cdot,dx)$ be a  transition kernel satisfying Hypothesis \ref{H:Kmeas}  for some $\alpha \in [0,\,1]$, and $H$ be a functional satisfying Hypothesis \ref{H:Hcont}. 
 $(Y, \P)$	fulfills the  martingale problem in Definition \ref{D:mtpb_hom_Mark} with respect to $C^2_b$, $\bar{\mathcal L}$ in \eqref{barL} and $y_0 \in \R$ if and only if, 
   for any $\tilde f \in C^2_b$, 
\begin{align}\label{mtg_pbY}
 &\tilde f(Y_{t}) - \tilde f(y_0) - \int_0^{\cdot} \bar L \tilde f(Y_{s}) ds - \int_0^{t}\sigma_0(Y_s) \bar H(s,Y^{-}) \normalcolor \tilde f'(Y_{s}) ds
\notag\\
&-   \int_{0}^t\int_{\R} (\tilde f(Y_{s-} + z) -\tilde f(Y_{s-})
-k(z)\,\tilde f'(Y_{s-})\,
)F(Y_{s-}, dz)ds
\end{align}
is an $(\mathcal F_t^Y)$-local martingale  under $\P$.
\end{remark}

For every $x \in \R$, we  define  $\mathcal H_x: w \mapsto h(x+w)-h(x)$ and its inverse function  $\mathcal H^{-1}_x: w \mapsto h^{-1}(h(x)+w)-x$.
\begin{remark}\label{R:5.18}
	$F(h(x), \cdot)$ is the push forward of $Q(x, \cdot)$ via $\mathcal H^{-1}_x$, so that $Q(x, \cdot)$ is the push forward of $F(h(x), \cdot)$ through $\mathcal H_x$.  
\end{remark}

\begin{theorem}\label{equiv_mtgpb_weaksol}
Let $\alpha \in [0,1]$.
Assume Hypotheses  \ref{H:h}, \ref{H:hBIS} and \ref{H:h2} with respect to $\alpha$.
 Let    $Q(\cdot,dx)$ be a  transition kernel satisfying
 Hypothesis \ref{H:Kmeas} with respect to  $\alpha$, and  $H$ be a functional verifying Hypothesis \ref{H:Hcont}.
Then
 $(X, \P)$	fulfills the  martingale problem in Definition \ref{D:mtpb_hom_Mark} with respect to ${\mathcal D}_{\mathcal L}$ in \eqref{barD}, $\mathcal L$ in \eqref{Ldistrib} and $x_0 \in \R $
if and only if $(Y = h(X), \P)$  fulfills the martingale problem in Definition \ref{D:mtpb_hom_Mark} with respect  to  $C^2_b$, $\bar {\mathcal{L}}$ in \eqref{barL}
 and $h(x_0)$.
 \end{theorem}

\proof
$(\Rightarrow)$  
Let $\tilde f \in C^2_b$ and set $f:= \tilde f \circ h$. Recalling that $h \in C^{1+\alpha}_{\textup{loc}}$, we have $f \in \mathcal D_{\mathcal L}$ by Proposition \ref{P:resold}-d). 
 $(X, \P)$	fulfills the  martingale problem in Definition \ref{D:mtpb_hom_Mark} with respect to ${\mathcal D}_{\mathcal L}$, $\mathcal L$ and $x_0$ if and only if, 
   for any $f \in \mathcal D_{{\mathcal L}}$,
\begin{align*}
f(X_{t}) - f(x_0) - \int_0^{t}(\mathcal{L} f)(s,X^{-}) ds
\end{align*}
is an $(\mathcal F_t^X)$-local martingale  under $\P$.
Setting  $y_0 = h^{-1}(x_0)$, this yields that 
$$
\tilde f(Y_t) - \tilde f(y_0)-\int_0^{t}
 ({\mathcal L}  f)(s,h^{-1}(Y^-))ds
$$
is an $(\mathcal F_t^X)$-local martingale  under $\P$, therefore also an $(\mathcal F_t^Y)$-local martingale, since $X$ and $Y$ have the same canonical filtration.
Using the form of $\mathcal L$ in \eqref{Ldistrib} and  Proposition \ref{P:resold}-d),  we get that
\begin{align}\label{martXPRIMA}
&\tilde f(Y_t) - \tilde f(y_0) - \int_0^t {L}^0 \tilde f(Y_s) ds - \int_0^{t}\sigma(Y_s) H(s,h^{-1}(Y^{-})) (h' \circ h^{-1})(Y_s)\tilde f'(Y_{s}) ds
\notag\\
&-  \int_0^t\int_{\R} [f(X_{s-} + w) - f(X_{s-})-k(w)\,f'(X_{s-})]Q(X_{s-}, \,dw)ds
\end{align}
is an $(\mathcal F_t^Y)$-local martingale  under $\P$. 
From  Remark \ref{R:5.18}, we have
\begin{equation*}
Q(x, A)=\int_{\R}\one_{A} \, (\mathcal H^{-1}_x(z))  F(h(x), dz). 
\end{equation*}
Therefore, we obtain  
\begin{align}
&\int_{\R}[f(X_{s-} + w)- f(X_{s-})-k(w)\,f'(X_{s-})] \,Q(X_{s-}, dw)\label{toplug}\notag\\
&=\int_{\R}[f(X_{s-} +\mathcal H^{-1}_{X_{s-}}(z))- f(X_{s-}) -k(\mathcal H^{-1}_{X_{s-}}(z))\,f'(X_{s-})] F(Y_{s-}, dz)\notag\\
&=\int_{\R}[f(h^{-1}(Y_{s-}) + h^{-1}(Y_{s-} + z)- h^{-1}(Y_{s-}))- f(h^{-1}(Y_{s-}))\notag\\
&\qquad  \qquad -k(h^{-1}(Y_{s-} + z)- h^{-1}(Y_{s-}))\,f'(h^{-1}(Y_{s-}))] \,F(Y_{s-}, dz)\notag\\
&=\int_{\R}[\tilde f(Y_{s-} + z)- \tilde f(Y_{s-}) -k(h^{-1}(Y_{s-} + z)- h^{-1}(Y_{s-}))\,\tilde f'(Y_{s-}) (h' \circ h^{-1})(Y_{s-})]  F(Y_{s-}, dz).	
\end{align}
Plugging \eqref{toplug} into \eqref{martXPRIMA} we get that
\begin{align}\label{martX}
&\tilde f(Y_t) - \tilde f(y_0) - \int_0^t {L}^0 \tilde f(Y_s) ds - \int_0^{t} \sigma(Y_s) H(s,h^{-1}(Y^{-})) (h' \circ h^{-1})(Y_s)\tilde f'(Y_{s}) ds
\notag\\
&-\int_0^t\int_{\R}[\tilde f(Y_{s-} + z)- \tilde f(Y_{s-}) -k(h^{-1}(Y_{s-}  + z)- h^{-1}(Y_{s-}))(h' \circ h^{-1})(Y_{s-})\,\tilde f'(Y_{s-}) ]  F(Y_{s-}, dz)ds
\end{align}
is an $(\mathcal F_t^Y)$-local martingale under $\P$.
Formula \eqref{martX} can be equivalently rewritten as
\begin{align*}
&\tilde f(Y_t) - \tilde f(y_0) - \int_0^t { L}^0 \tilde f(Y_s) ds 		- \int_0^{t} \sigma_0(Y_s) 
H(s,h^{-1}(Y^{-}))
 \tilde f'(Y_{s}) ds\\
&-  \int_0^t\int_{\R} [\tilde f(Y_{s-} + z) -\tilde f(Y_{s-})-k(z)\tilde f'(Y_{s-})]  F(Y_{s-}, dz)\,ds
\\
&-  \int_0^t\tilde f'(Y_{s-}) (h' \circ h^{-1})(Y_{s-})\int_{\R} [k(z)\,(h^{-1})'(Y_{s-})-k(h^{-1}(Y_{s-} + z)- h^{-1}(Y_{s-}))]F(Y_{s-}, dz)ds,
\end{align*}
which provides formula \eqref{mtg_pbY}
  with the operators $ {\bar L}$  and $\bar H$ given respectively  by \eqref{tildeL} an \eqref{barH}.
This finally shows that $(Y, \P)$  fulfills the martingale problem in Definition \ref{D:mtpb_hom_Mark} related to  $C^2_b$, $\bar {\mathcal{L}}$ in \eqref{barL}
 and $h(x_0)$.

\medskip

\noindent $(\Leftarrow)$ 
%
Let $f \in {\mathcal D}_{\mathcal L}$
    and  set $\phi = f \circ h^{-1}$. By Proposition \ref{P:resold}-d)  $\phi \in C^2_b$. 
Then, by assumption, 
\begin{align*}
&\phi(Y_t) - \phi(h(x_0)) 
- \int_0^{t} L^0 \phi(Y_{s}) ds - \int_0^{t} \sigma(Y_s) H(s,h^{-1}(Y^{-})) (h' \circ h^{-1})(Y_s)\phi'(Y_{s}) ds		\\
&-  \int_0^{t} \int_{\R}(\phi(Y_{s-} + z) -\phi(Y_{s-})
-k(z)\,\phi'(Y_{s-}))F(Y_{s}, dz) ds\\
&+\int_0^{t}\int_{\R} \phi'(Y_{s}) (h' \circ h^{-1})(Y_{s})[k(h^{-1}(Y_{s} + z)- h^{-1}(Y_{s}))-(h^{-1})'(Y_{s})\,k(z)]F(Y_{s}, dz) ds 
\end{align*}
is an $(\mathcal F_t^Y)$-local martingale under $\P$, that in turn gives that
\begin{align}\label{martXPRIMABIS}
& f(X_t) - f(X_0) - \int_0^t {L} f(X_s) ds - \int_0^{t} \sigma(X_s) H(s,X^{-})  f'(X_{s}) ds
\notag
\\
&-  \int_0^t \int_{\R} [\phi(Y_{s-} + z) - \phi(Y_{s-})-\phi'(Y_{s-}) (h' \circ h^{-1})(Y_{s-})k(h^{-1}(Y_{s-} + z)- h^{-1}(Y_{s-}))]\,F(Y_{s-}, dz)ds \end{align}
is an $(\mathcal F_t^X)$-local martingale under $\P$. 
At this point, using \eqref{identity_nu_nuY_F},  we get 
\begin{align}
&\int_{\R}[\phi(Y_{s-} + z)- \phi(Y_{s-}) -\phi'(Y_{s-}) (h' \circ h^{-1})(Y_{s-})\,k(h^{-1}(Y_{s-} + z)- h^{-1}(Y_{s-}))] \, F(Y_{s-}, dz)\notag\\
&=\int_{\R}[f(X_{s-} + w)- f(X_{s-})-f'(X_{s-})\,k(w)] \,Q(X_{s-}, dw).\label{toplugBIS}
\end{align}
Plugging \eqref{toplugBIS} into \eqref{martXPRIMABIS} we get the result.
 \qed
 
 \subsection{Weak Dirichlet property}
 
 The notion of characteristics of weak Dirichlet processes was introduced in   Section 3.3 in \cite{BandiniRusso_RevisedWeakDir}, extending the classical one for semimartingales, see Appendix \ref{A:sem}.
  We will denote by $X^c$ the unique continuous local martingale component of $X$, see Proposition 3.2 in \cite{BandiniRusso_RevisedWeakDir}.

Below,  $\check Y$ replaces $\check X$ in the role of canonical process.
  \begin{proposition}\label{R:Yspecsem}
 Let $\alpha \in [0,1]$.
Assume Hypotheses  \ref{H:h}, \ref{H:hBIS} and \ref{H:h2} with respect to $\alpha$.
Let    $Q(\cdot,dx)$ be a  transition kernel satisfying Hypotheses  \ref{H:Kmeas} 
  with respect to $\alpha$, and  $H$ be a functional satisfying  Hypothesis \ref{H:Hcont}. 
If
$(X, \P)$ is  a  solution to the martingale problem in Definition \ref{D:mtpb_hom_Mark} related to $\mathcal {D}_{\mathcal L}$ in \eqref{barD}, $\mathcal L$ in \eqref{Ldistrib} and $x_0 \in \R$,  then the following holds.
    \begin{enumerate}
    	\item 
$Y=h(X)$ is a  semimartingale with characteristics 
$B = \int_0^\cdot (b(\check Y_s) + \sigma_0(\check Y_s)  \bar H(s,\check Y^{-})) ds$, $C = \int_0^\cdot c(\check Y_s) ds$, $\tilde \nu(ds\,dz) = F(\check Y_{s}, dz)ds$, where $\sigma_0$,  $b$ and $\bar H$  are defined respectively in \eqref{sigma0},   \eqref{b(y)} and \eqref{barH}, $F(y, dz)$ is the measure introduced in  \eqref{identity_nu_nuY_F}, and 
			$
			c(y) := \sigma_0^2(y)
			$.
			\item 
 	  $X$ is a  weak Dirichlet process of finite quadratic variation with characteristic $\nu(ds\, dw) = Q(\check X_{s-}, dw) ds$. 
 	  \item  $\langle X^c, X^c\rangle = \int_0^\cdot\sigma^2(X_s) ds $.
 			    \end{enumerate} 
  \end{proposition} 

\proof 
1.  It is a direct consequence of Theorems \ref{equiv_mtgpb_weaksol} and \ref{T: equiv_mtgpb_semimart}.



 2. 
By definition $X = h^{-1}(Y)$, with $h^{-1} \in C^1$. 
By item 1., $Y$ is a semimartingale, so it is  a  weak Dirichlet process  of finite quadratic variation.
 In particular, $X$ has finite quadratic variation, see
Lemma \ref{L:app1}-1.
Moreover, 
we can  apply   Theorem  3.36 in \cite{BandiniRusso_RevisedWeakDir}
to  $h^{-1}(Y)$, so that 
$X$ turns out to be a  weak Dirichlet process.
Finally, by item 1. and \eqref{identity_nu_nuY_F}, 
\begin{align*}
	\tilde \nu(ds,A) &= F(h(\check X_{s-}), A)ds\\
	&= \int_{\R}\one_{A} \, (h(h^{-1}(h(\check X_{s-}))+ w)- h(h^{-1}(h(\check X_{s-})))) \, Q(h^{-1}(h(\check X_{s-})), dw)ds\\
	&=\int_{\R}\one_{A} \, (h(\check X_{s-}+ w)- h(\check X_{s-})) \, Q(\check X_{s-}, dw)ds.
\end{align*}
Then, by Remark 3.41   in \cite{BandiniRusso_RevisedWeakDir} with $v(t,y)= h^{-1}(y)$, $\nu^Y= \nu$ and  $\nu^X =\tilde  \nu$,
the characteristic $\nu$ of $X$ is given by 
  \begin{align*}
  \nu(A, ds)&= \int_\R \one_{A} \, (h^{-1}(h(\check X_{s-})+z)-\check X_{s-}) \, \tilde \nu(ds,\,dz) 	\\
  &=\int_\R \one_{A} \, (h^{-1}(h(\check X_{s-})+h(\check X_{s-}+ w)- h(\check X_{s-}))-\check X_{s-}) \, Q(\check X_{s-}, dw)ds\\
  &=\int_\R \one_{A} \, (h^{-1}(h(\check X_{s-}+ w))-\check X_{s-}) \, Q(\check X_{s-}, dw)ds\\
  &=\int_\R \one_{A}(w) \, Q(\check X_{s-}, dw)ds.
  \end{align*}
  
3. From item 1, 
$$
C \circ Y = \int_0^\cdot (\sigma^2 h')(h^{-1}(h( X_s))) ds=\int_0^\cdot |h'(X_s)|^2\sigma^2(X_s) ds.
$$
On the other hand, by formula (3.45) in  Remark 3.42-(i) in \cite{BandiniRusso_RevisedWeakDir}, 
$$
C \circ Y 
= \int_0^\cdot |h'(X_s)|^2   d \langle  X^c, X^c\rangle_s,
$$
and the conclusion follows.
\endproof

\begin{remark}\label{R:3.14}
If
$(X, \P)$ is  a  solution to the martingale problem in Definition \ref{D:mtpb_hom_Mark} related to $\mathcal {D}_{\mathcal L}$  in \eqref{barD}, $\mathcal L$ in \eqref{Ldistrib},  and $x_0\in \R$, then it is not necessarily a  Dirichlet process. 

Consider for instance the case $X = W + S$ with $W$ a Brownian motion and $S$ an $\gamma$-stable L\'evy process with $\gamma \in (0,1)$.
This can be seen as a trivial solution of our martingale problem with $\sigma\equiv 1$ and $Q(y, dx)= Q_0(dx) =  |x|^{-1-\gamma} dx$. We remark that $X$ is a Dirichlet process if and only if $S$ is a Dirichlet process.
Assume ab absurdo that  $S$ is a Dirichlet process. Since $S$ is also a semimartingale,  then 
 $S$ is  special 
semimartingale, see Lemma \ref{L:3.2} and Proposition 5.14 in \cite{BandiniRusso1}.
However,  $x \one_{|x| >1} \star Q_0 = + \infty$, and therefore it cannot be a special semimartingale,  see 
Proposition 2.29, Chapter II, in \cite{JacodBook}.
Notice that, in the case $\gamma \in [1,2)$, $S$ instead  is a special semimartingale because $x \one_{|x| >1} \star Q_0 < + \infty$. 

We will state and prove new results  on Dirichlet processes in  Section \ref{S:Dirichlet}. 
\end{remark}

\section{Well-posedness of the martingale problem} \label{S5}
In order to formulate the well-posedness of the martingale problem we will make use of the following hypothesis about some transition kernel $Q(\cdot, dx)$. 
\begin{hypothesis}\label{H:totalvar}
For some $\alpha \in [0,\,1]$,
\begin{align}\label{measure}
y \mapsto  (1 \wedge |x|^{1+\alpha}) \, Q(y,dx)\quad \textup{is  continuous in the total variation topology}.
\end{align}
\end{hypothesis}

\begin{remark}
	According to Remark \ref{R:Levystable}, Hypothesis \ref{H:totalvar} is trivially verified in  the case of $Q(y,dx) = Q_0(dx) = |x|^{-1-\gamma} dx$ if $\alpha > \gamma-1$,  being the measure-valued function \eqref{measure} constant.
\end{remark}

\begin{remark}\label{R:app}
\begin{itemize}
	\item [(i)]If  
Hypothesis  \ref{H:totalvar} holds true for some $\alpha \in [0,\,1]$, then  
$
y \mapsto  (1 \wedge |x|^{2}) \, Q(y,dx)
$
is continuous in the total variation topology.
\item [(ii)] Item (i)  in turn  implies that 
$
y \mapsto \int_{B} (1 \wedge |x|^{2}) \, Q(y,dx)
$
is continuous 
 for all $B\in \mathcal B(\R)$.
\end{itemize}
\end{remark}
 
 We consider again the functions $\Sigma$ and $h$  introduced respectively in \eqref{Sigma} and  in Proposition  \ref{P:equiv}.
  We will make the following additional assumption. 
\begin{hypothesis}\label{H:Sigma}
	$\Sigma$ is bounded and is $\alpha$-H\"older continuous  in the whole space for some  $\alpha \in [0,\,1]$ (where $0$-H\"older continuous means uniformly continuous).
\end{hypothesis}
\begin{remark}\label{R:4.5}
\begin{itemize}
	\item [(i)] Under Hypothesis \ref{H:Sigma}, $h'$ is upper and lower bounded as well.
	\item [(ii)] Hypothesis \ref{H:Sigma} implies Hypotheses \ref{H:hBIS} and  \ref{H:h2}.
	\item[(iii)] For some  $\alpha \in (0,\,1)$, 
 Hypothesis \ref{H:Sigma} is equivalent to ask that $\Sigma$ belongs to the Besov space $\mathcal C^\alpha$, see e.g. Section 2.7 in \cite{bahouri}.
\end{itemize}	
\end{remark}
\normalcolor
We start by considering the Markovian case. 

\begin{proposition}
\label{P47}
Let $\alpha \in [0,1]$.
Let $L$ be an operator of the form \eqref{Lbeta}  with $\sigma$ bounded.
Assume Hypotheses  \ref{H:h} and \ref{H:Sigma} with respect to $\alpha$.
Let   $Q(\cdot,dx)$ be a  transition kernel satisfying Hypotheses \ref{H:Kmeas}, and \ref{H:totalvar} with respect to $\alpha$.
Then  the martingale problem in Definition \ref{D:mtpb_hom_Mark}
related to   $\mathcal D_{\mathcal L}$ in \eqref{barD},  $\mathcal L$ in \eqref{Ldistrib}  with $H \equiv 0$ and $x_0 \in \R$ admits existence and uniqueness.
\end{proposition} 
\proof
By Theorem \ref{equiv_mtgpb_weaksol}, existence and uniqueness of the  martingale problem in Definition \ref{D:mtpb_hom_Mark} with respect to ${\mathcal D}_{\mathcal L}$, $\mathcal L$ with $H \equiv 0$ and $x_0$
is equivalent to existence and uniqueness of the martingale problem in Definition \ref{D:mtpb_hom_Mark} related to $C^2_b$
$\bar {\mathcal{L}}$
in \eqref{barL} with $\bar H\equiv 0$
  and $h(x_0)$. 
On the other hand, by Theorem \ref{T: equiv_mtgpb_semimart}, 
 $Y=h(X)$ is a solution to the latter martingale problem if and only if it a 
 semimartingale with local characteristics 
$B= \int_0^\cdot b(\check Y_s) ds$, $ C = \int_0^\cdot c(\check Y_s) ds$, $\tilde \nu(ds\,dz) = F(\check Y_{s}, dz)ds$, with, for every $y \in \R$,  
\begin{align*}
F(y, A)&:=\int_{\R }\one_{A} \, (h(h^{-1}(y)+ w)- h(h^{-1}(y))) \, Q(h^{-1}(y), dw), \quad A \subseteq \R,\\
b(y) &:= (h' \circ h^{-1})(y)\int_{\R}[(h^{-1})'(y)\,k(z)-k(h^{-1}(y + z)- h^{-1}(y))]\,F(y, dz),\\
c(y) &:= (\sigma h')^2(h^{-1}(y)).
\end{align*}
The result will then follow by using 	Theorem \ref{T:2.34}, provided 
we  verify Hypothesis \ref{H_exis_uni_mtgpb}  for  $b$, $c$  and $F(\cdot,dz)$, i.e. that  
\begin{itemize}
\item [(i)]	$b$ is bounded;
\item  [(ii)]	$c$ is  bounded, continuous, and not vanishing at zero;
\item [(iii)] the function $y \mapsto \int_{B} (1\wedge |z|^2) F(y,dz)$ is bounded and continuous for all $B \in  \mathcal B(\R)$. 
\end{itemize}

We start by item (ii). 
Recall that by Remark \ref{Rh}  and Proposition \ref{P:equiv}, we can take $h \in C^{1}$,   $h'$ being bounded and $h^{-1}$ being continuous.
Since $\sigma$ is continuous, this implies that  the function  $c$   is continuous as well.  Moreover, 
since $\sigma$ is bounded, $c$ is also bounded. Finally, $c$ is not vanishing at zero by formula \eqref{h'} and the fact that $\sigma$ is never zero.

We then prove that 
\begin{itemize}
\item [(iii)'] the function 
$y \mapsto (1 \wedge |z|^{1+\alpha}) F(y,dz)$ is bounded and continuous in the total variation norm.
\end{itemize}
In particular, this would imply item (iii), see Remark 
\ref{R:app}.
 We have 
  \begin{align*}
 (1 \wedge |z|^{1+\alpha})\,F(y, dz)
&=
 (1 \wedge |h(h^{-1}(y) + w) -h(h^{-1}(y))|^{1+\alpha})\,Q(h^{-1}(y), dw) \\
&=  
 (1 \wedge (\psi(y, w)\,|w|)^{1+\alpha}) \,Q(h^{-1}(y), dw):= I(y;dw)
\end{align*}
with 
\begin{equation*}
\psi(y, w):=\int_0^1 h'(h^{-1}(y) + a w) da.	
\end{equation*}
Since $h'$ is bounded, there is  a constant $C_1$ such that  $\psi^{1+\alpha} \leq C_1$.

Let us first  prove the boundedness of the map $y \mapsto I(y;dw)$. We have $I(y;dw) = I_1(y;dw) + I_2(y;dw)$ with
\begin{align*}
	I_1(y;dw) &:= \one_{\{0 < |w|^{1+\alpha} \leq \frac{1}{C_1}\}} (1 \wedge (\psi(y, w)\,|w|)^{1+\alpha}) \,Q(h^{-1}(y), dw),\\
	I_2(y;dw) &:= \one_{\{|w|^{1+\alpha} > \frac{1}{C_1}\}} (1 \wedge (\psi(y, w)\,|w|)^{1+\alpha}) \,Q(h^{-1}(y), dw).
\end{align*}
For  $\ell_1(w) :=\one_{\{0 < |w|^{1+\alpha} \leq \frac{1}{C_1}\}}$ and $\ell_2(w) :=\one_{\{|w|^{1+\alpha} > \frac{1}{C_1}\}}$, we set
\begin{align*}
 	\tilde Q^{\ell_1}(h^{-1}(y), dw)&:=\one_{\{0 < |w|^{1+\alpha} \leq \frac{1}{C_1}\}} |w|^{1+\alpha} Q(h^{-1}(y), dw),\\
 	\tilde Q^{\ell_2}(h^{-1}(y), dw)&:=\one_{\{|w|^{1+\alpha} > \frac{1}{C_1}\}} Q(h^{-1}(y), dw).
\end{align*}
For every $y \in \R$, 
\begin{align*}
	||I_1(y;dw)||_{var}&\leq C_1 \int_\R \one_{\{0 < |w|^{1+\alpha} \leq \frac{1}{C_1}\}} |w|^{1+\alpha}  \,Q(h^{-1}(y), dw)= C_1 \sup_{z \in \R}||\tilde Q^{\ell_1}(z, dw)||_{var}
\end{align*}
and 
$$
	||I_2(y;dw)||_{var} 
\leq \sup_{z \in \R} ||\tilde Q^{\ell_2}(z, dw)||_{var},
$$
whereas previous supremum are finite  by Lemma \ref{L:techResult}. 

Let us now prove the continuity of the map $y \mapsto I(y;dw)$. Let $(y_n)$ be a real sequence converging to $y_0 \in \R$.
 We have 
$$
I(y_n;dw)-I(y_0;dw) = J_1(y_n, y_0;dw) + J_2(y_n, y_0;dw)
$$
with 
\begin{align*}
J_1(y_n, y_0;dw)&:= 
 (1 \wedge (\psi(y_n, w)\,|w|)^{1+\alpha})  \,[Q(h^{-1}(y_n), dw)-Q(h^{-1}(y_0), dw)],\\
J_2(y_n, y_0;dw)&:=
\{(1 \wedge (\psi(y_n, w)\,|w|)^{1+\alpha}) \,-(1 \wedge (\psi(y_0, w)\,|w|)^{1+\alpha})\} \,Q(h^{-1}(y_0), dw).
\end{align*}
Concerning $J_1$, we have $J_1= J_1'+ J_1''$, where
\begin{align*}
J_1'(y_n, y_0;dw)&:= \one_{\{0 < |w|^{1+\alpha} \leq \frac{1}{C_1}\}} (1 \wedge (\psi(y_n, w)\,|w|)^{1+\alpha})  \,[Q(h^{-1}(y_n), dw)-Q(h^{-1}(y_0), dw)],\\
J_1''(y_n, y_0;dw)&:= \one_{\{|w|^{1+\alpha} > \frac{1}{C_1}\}} (1 \wedge (\psi(y_n, w)\,|w|)^{1+\alpha})  \,[Q(h^{-1}(y_n), dw)-Q(h^{-1}(y_0), dw)].
\end{align*}
 We get 
\begin{align*}
	J_1'(y_n, y_0;dw)&= \one_{\{0 < |w|^{1+\alpha} \leq \frac{1}{C_1}\}} \psi^{1+\alpha}(y_n, w)\,|w|^{1+\alpha}  \,[Q(h^{-1}(y_n), dw)-Q(h^{-1}(y_0), dw)]
	\end{align*}
	so that 
	\begin{align*}
	||J_1'(y_n, y_0;dw)||_{var}
	&\leq C_1 ||\tilde Q^{\ell_1}(h^{-1}(y_n), dw) - \tilde Q^{\ell_1}(h^{-1}(y_0), dw)||_{var}
\end{align*}
and analogously 
 \begin{align*}
	||J_1''(y_n, y_0;dw)||_{var}&\leq   ||\one_{\{|w|^{1+\alpha} > \frac{1}{C_1}\}} [Q(h^{-1}(y_n), dw)-Q(h^{-1}(y_0), dw)]||_{var}\\
	&\leq ||\tilde Q^{\ell_2}(h^{-1}(y_n), dw) - \tilde Q^{\ell_2}(h^{-1}(y_0), dw)||_{var}.
\end{align*}
The convergence of both terms follows by Lemma \ref{L:techResult} applied respectively to $\tilde Q^{\ell_1}(h^{-1}(y_0), dw)$ and $\tilde Q^{\ell_2}(h^{-1}(y_0), dw)$, and taking into account  the continuity of $h^{-1}$.

Regarding $J_2$ 
we have $J_2= J_2'+ J_2''$
with
\begin{align*}
J_2'(y_n, y_0;dw):=
& \{(1 \wedge (\psi(y_n, w)\,|w|)^{1+\alpha}) \,-(1\wedge (\psi(y_0, w)\,|w|)^{1+\alpha})\} \one_{\{0 < |w|^{1+\alpha} \leq \frac{1}{C_1}\}}(w)\,Q(h^{-1}(y_0), dw),\\
J_2''(y_n, y_0;dw):= 
&\{(1 \wedge (\psi(y_n, w)\,|w|)^{1+\alpha}) \,-(1 \wedge (\psi(y_0, w)\,|w|)^{1+\alpha})\} \one_{\{|w|^{1+\alpha} > \frac{1}{C_1}\}}(w)\,Q(h^{-1}(y_0), dw).
\end{align*}
Notice that 
\begin{align*}
J_2'(y_n, y_0;dw)= (\psi(y_n, w)\,  \,-\psi(y_0, w)) \one_{\{0 < |w|^{1+\alpha} \leq \frac{1}{C_1}\}}(w)|w|^{1+\alpha}\,Q(h^{-1}(y_0), dw),
\end{align*}
so that 
\begin{align*}
||J_2'(y_n, y_0;dw)||_{var}\leq  \int_\R |\psi(y_n, w)\,  \,-\psi(y_0, w)| \one_{\{0 < |w|^{1+\alpha} \leq \frac{1}{C_1}\}}(w)|w|^{1+\alpha}\,Q(h^{-1}(y_0), dw).
\end{align*}
On the other hand,
\begin{align*} 
||J_2''(y_n, y_0;dw)||_{var}\leq \int_\R \{(1 \wedge (\psi(y_n, w)\,|w|)^{1+\alpha}) \,-(1 \wedge (\psi(y_0, w)\,|w|)^{1+\alpha})\} \one_{\{|w|^{1+\alpha} > \frac{1}{C_1}\}}(w)\,Q(h^{-1}(y_0), dw).
\end{align*}
Therefore
$||J_2'(y_n, y_0;dw)||_{var}$ and $||J_2''(y_n, y_0;dw)||_{var}$
 converge to zero by the Lebesgue dominated convergence theorem, taking into account respectively the finiteness  $\tilde Q^{\ell_1}(h^{-1}(y_0), dw)$ and $\tilde Q^{\ell_2}(h^{-1}(y_0), dw)$ due to Lemma \ref{L:techResult}, and the continuity of $h', h^{-1}$.
This proves (iii)'.


Finally, let us prove item (i). We first notice that 
 \begin{align}\label{to_specify}
b(y)&=(h' \circ h^{-1})(y)\int_{\R}[(h^{-1})'(y)\,k(z)-k(h^{-1}(y + z)- h^{-1}(y))]\,F(y, dz)\notag\\
&=(h' \circ h^{-1})(y)\int_{\R}[(h^{-1})'(y) \,k(z)-k(z\,\bar \psi(y,z))]\,F(y, dz),
 \end{align}
 with 
$
\bar \psi(y, z):=\int_0^1 (h^{-1})'(y + a\,z)\,da	
$. 
Also in this case we can find a  constant $\bar C_1 \geq 1$ such that $\psi \leq \bar C_1$.
For some $R \in (0,\,1)$, define   $\mathcal B_R:=\{z \in \R:\,\, |z| \leq R\}$ as  the neighborhood  of $z=0$  on which $k(z) =z$. We also introduce $\bar {\mathcal B}:=\{z \in \R:\,\, |z| \leq \frac{R}{\bar C_1}\}\subset \mathcal B_R$. 
Identity \eqref{to_specify} reads 
\begin{align}\label{expansionsbk}
	b(y)&= (h' \circ h^{-1})(y)\int_{\R}\Big[\int_0^1 ((h^{-1})'(y)-(h^{-1})'(y + a\,z) )\,da\Big] z\,\, \one_{\bar {\mathcal B}}(z)\,F(y, dz)\notag\\
	&+(h' \circ h^{-1})(y)\int_{\R}[(h^{-1})'(y) \,k(z)-k(z\,\bar \psi(y, z)) ]\,\one_{\bar {\mathcal B}^c}(z)\,F(y, dz).
	\end{align}

In the sequel we suppose  $\alpha \in (0,1]$, the case $\alpha=0$ needs some easy adaptation.
Concerning the boundedness of $b$, we first notice that by \eqref{h'} together with  Hypothesis \ref{H:Sigma}, for every $a \in [0,1]$, 
\begin{align*}
	|(h^{-1})'(y)-(h^{-1})'(y + a\,z)| \,\one_{\bar {\mathcal B}}(z) &= |e^{\Sigma(h^{-1}(y))} -e^{\Sigma(h^{-1}(y+az))}| \one_{\bar {\mathcal B}}(z) \\
	&\leq C_2 \, e^{||\Sigma||_\infty} |h^{-1}(y)-h^{-1}(y+az)|^\alpha \one_{\bar {\mathcal B}}(z)\\
	&\leq C_2 \, e^{(1+\alpha||\Sigma||_\infty)} |z|^\alpha \one_{\bar {\mathcal B}}(z), 
\end{align*}
where $C_2$ is a  H\"older constant for $\Sigma$.
Therefore  by \eqref{expansionsbk}
 \begin{align*}
|b(y)|
&\leq  ||h' ||_{\infty} 
C_2 \, e^{(1+\alpha||\Sigma||_\infty)} \int_{\R }  |z|^{1+\alpha} \, \one_{\bar {\mathcal B}}(z)\, F(y, dz)+  ||h'||_{\infty}  ||k||_{\infty}(1+ ||(h^{-1})'||_{\infty}) \int_{\R}\one_{\bar {\mathcal B}^c}(z)\,F(y, dz),
 \end{align*}
and the conclusion follows   by Lemma \ref{L:techResult} applied to $\ell_1(z) = \one_{\mathcal B}(z)$ and $\ell_2(z) = \one_{\mathcal B^c}(z)$. 
\endproof

We finally can state the general existence and uniqueness theorem for the possibly path-dependent case.

\begin{theorem} \label{T47}
Let $\alpha \in [0,1]$.
Let  $L$ be an operator of the form \eqref{Lf}   with $\sigma$ bounded. Assume Hypotheses  \ref{H:h} and \ref{H:Sigma} with respect to $\alpha$.
Let   $Q(\cdot,dx)$ be a  transition kernel satisfying Hypotheses \ref{H:Kmeas} and \ref{H:totalvar}
with respect to  $\alpha$,  and $H$ be a functional satisfying Hypothesis \ref{H:Hcont}. 
Then 
existence and uniqueness holds for the martingale problem in Definition \ref{D:mtpb_hom_Mark}
related to  $\mathcal D_{\mathcal L}$ in \eqref{barD}, $\mathcal L$ in \eqref{Ldistrib}  and  $x_0 \in \R$.
\end{theorem}
\proof
\emph{Step 1.}
%
Let $X$ be an $(\mathcal F^X_t)$-weak Dirichlet 
with characteristic $\nu(ds\,dx)$ such that  $(1 \wedge |x|^2)\star (\nu \circ X) \in \mathcal A_{\textup{loc}}$, 
and with $(\mathcal F^X_t)$-continuous local martingale $X^{c}$ under $\P$ such that $\langle X^{c}, X^{c}\rangle=\int_0^\cdot \sigma^2(X_s) ds$.
  We  set 
\begin{equation*}
W_t :=\int_0^t \frac{1}{\sigma(X_s)} d X^c_s, \quad t \in [0,T]. 	
\end{equation*}
Consequently $W$ is an $(\mathcal F^X_t)$-local martingale with $\langle W, W\rangle_t=t$, and therefore by L\'evy's characterization theorem, $W$ is an $(\mathcal F^X_t)$-Brownian motion. 
Let $H$ be a functional defined on $D_{-}(0,T)$
 satisfying Hypotheses \ref{H:Hcont}. We define 
\begin{equation}\label{tildeW}
	\tilde W_t := W_t - \int_0^t  H(s,X^{-}) ds, \quad t \in [0,T]. 
\end{equation}
Then, by the Novikov condition, 
\begin{align*}
\kappa_t :=\exp{\Big\{\int_0^t  H(s,X^{-}) dW_s- \frac{1}{2}\int_0^t| H(s,X^{-})|^2 ds\Big \}}, \quad t \in [0,T], 	
\end{align*}
is an $(\mathcal F^X_t)$-martingale. 
By Girsanov's theorem, $\tilde W$ is an $(\mathcal F^X_t)$-Brownian motion under the probability $\Q$ defined by 
\begin{equation}\label{Q}
	d \Q= \kappa_T d \P.
\end{equation}
Let $f \in \mathcal D_{\mathcal L}$, and 
set $\eta_s(x):=f(X_{s-} +x)-f(X_{s-})$   and 
$$
\xi_s(x) := \eta_s(x) \star (\mu^X-(\nu \circ X)).
$$
The process $\xi$ is an $(\mathcal F^X_t)$-purely discontinuous
local martingale under $\P$, see considerations in Definition 1.27-(ii), Chapter II,  in \cite{JacodBook}.
In particular $\langle \xi, M \rangle = 0$  for every continuous local 
martingale $M$.
We claim  that 
\begin{equation}\label{claim}
\text{$\xi$    remains an $(\mathcal F^X_t)$-local martingale under $\Q$.}
\end{equation}
Indeed, set $\tau_n :=\inf\{t \in [0,T]:\,\,|X_{t-}|>n\}$.
We recall that the c\`agl\`ad process $(X_{t-})$ is locally bounded.
Then the process $\xi^n:= \xi \one_{[0, \tau_n]}$ is a (square integrable) martingale under $\P$. 
As a matter of fact,  $\eta_s(x)\one_{[0,\tau_n]}(s)\in \mathcal L^2(\mu^X)$ (and in particular belongs to $\mathcal G^2(\mu^X)$, see the end of Section 2 in \cite{BandiniRusso_RevisedWeakDir}) since  
\begin{align*}
\eta^2_s(x) \one_{\{|x|>1\}}&\leq 4||f||^2_\infty,\\
\eta^2_s(x) \one_{[0,\tau_n]}(s)\one_{\{|x|\leq 1\}}  &\leq ||f'(\cdot)\one_{[-(n+1), n+1]}(\cdot)||_\infty\, x^2 \one_{\{|x|\leq 1\}}.
\end{align*}
To prove that $\xi^n$ remains an $(\mathcal F^X_t)$-martingale under $\Q$, we need to show that, for every    $\mathcal F^X_s$-measurable random variable $F$,
$
\E^\Q[(\xi^n_t-\xi^n_s)F]=0.
$
Indeed, the left-hand side gives 
\begin{align*}
\E^\P[\kappa_T (\xi^n_t-\xi^n_s)F]= \E^\P[(\kappa_t-\kappa_s) (\xi^n_t-\xi^n_s)F]=\E^\P[(\langle \kappa, \xi\rangle^n_t -\langle \kappa, \xi^n \rangle_s) F]=0,
\end{align*}
since $\langle \kappa, \xi^n\rangle=0$, being $\xi^n$ an $(\mathcal F^X_t)$-purely discontinuous local martingale.
This shows that $\xi$ is an $(\mathcal F^X_t)$-local martingale under $\Q$.

\emph{Step 2: existence.}
 Let $(X, \P)$ be a solution to the martingale problem in Definition \ref{D:mtpb_hom_Mark}
related to   $\mathcal D_{\mathcal L}$ in \eqref{barD}, $\mathcal L$ in \eqref{Ldistrib} with $H \equiv 0$ and  $x_0 \in \R$. 
By Proposition \ref{R:Yspecsem}  with $H \equiv 0$, $X$
is  an $(\mathcal F^X_t)$-weak Dirichlet 
with characteristic $\nu(ds\,dx)= Q(\check X_{s-},dx)ds$, 
and with $(\mathcal F^X_t)$-continuous local martingale $X^{c}$ under $\P$ such that $\langle X^{c}, X^{c}\rangle=\int_0^\cdot \sigma^2(X_s) ds$.
By the uniqueness of the decomposition for special weak Dirichlet processes and Corollary 3.37 in \cite{BandiniRusso_RevisedWeakDir}, for every $f \in \mathcal D_{\mathcal L}$ we have  
\begin{align}\label{martpb_new}
&f(X_{\cdot}) - f(x_0) - \int_0^{\cdot}(\mathcal{L} f)(s,X^{-}) ds \notag\\
&= \int_0^\cdot (f'\sigma)(X_s) dW_s + \int_0^\cdot (f(X_{s-} +x)-f(X_{s-})) (\mu^X(ds\,dx)-Q(X_{s-},dx)ds).
\end{align}
Plugging  in \eqref{martpb_new} the process $\tilde W$ defined in \eqref{tildeW}, we get 
\begin{align*}
&f(X_{\cdot}) - f(x_0) - \int_0^{\cdot} (\mathcal{L} f)(s,X^{-}) ds -\int_0^\cdot (f'\sigma)(X_s) H(s,X^{-}) ds \notag\\
&= \int_0^\cdot (f'\sigma)(X_s) d\tilde W_s + \int_0^\cdot (f(X_{s-} +x)-f(X_{s-})) (\mu^X-\nu\circ X)(ds\,dx).
\end{align*} 
Let $\Q$ be the probability constructed in \eqref{Q}.  By \eqref{claim} in Step 1.  
$$\int_0^\cdot (f'\sigma)(X_s) d\tilde W_s + \int_0^\cdot (f(X_{s-} +x)-f(X_{s-})) (\mu^X-\nu\circ X)(ds\,dx)
$$
is an $(\mathcal F^X_t)$-local martingale under $\Q$.
Therefore,  $(X, \Q)$ is proved to be a solution to the martingale problem in the statement.

\emph{Step 3: uniqueness.}
Let $(X^i, \P^i)$, $i=1,2$, be two solutions of the martingale problem in Definition \ref{D:mtpb_hom_Mark}  related to   $x_0 \in \R$,  $\mathcal D_{\mathcal L}$ in \eqref{barD}, and  $\mathcal L$ in \eqref{Ldistrib}.
By Proposition \ref{R:Yspecsem}-(2)(3), $X^i$ is an $(\mathcal F^{X^i}_t)$-weak Dirichlet with characteristic $\nu(ds\,dx)= Q(\check X_{s-},dx)ds$, with $(\mathcal F^{X^i}_t)$-local martingale $X^{i,c}$ under $\P^i$ such that $\langle X^{i,c}, X^{i,c}\rangle=\int_0^\cdot \sigma^2(X^i_s) ds$.
Consequently, by L\'evy's characterization theorem, 
$$
W^i:= \int_0^\cdot \frac{1}{\sigma^2(X_s^i)}d X_s^{i,c}, \quad t \in[0,T],
$$
 is an $(\mathcal F^{X^i}_t)$-Brownian motion. 
We define the $\P^i$-martingale 
$$
\kappa_t^i := {\Big\{-\int_0^t H(s,X^{i-}) dW_s- \frac{1}{2}\int_0^t|H(s,X^{i-})|^2 ds\Big \}}, \quad t \in [0,T], 
$$
and the probability $\Q^i$ such that $d \Q^i= \kappa_T^i d \P^i$. By Girsanov's theorem, under $\Q^i$, 
$$
B_t^i := W_t^i + \int_0^t H(s,X^{i-}) ds 
$$
is a Brownian motion.
By formula \eqref{claim} in  Step 1 (replacing $H$ with $-H$),   $(f(X^i_{s-} +\cdot)-f(X^i_{s-})) \star (\mu^{X^i}-\nu\circ {X^i})$ remains an $(\mathcal F^{X^i}_t)$-martingale under $\Q^i$.

 Therefore, $(X^i, \Q^i)$ solves the martingale problem in Definition \ref{D:mtpb_hom_Mark}  related to    $\mathcal D_{\mathcal L}$ in \eqref{barD},   $\mathcal L$ in \eqref{Ldistrib} with $H \equiv 0$ and  $x_0 \in \R$.

 By the uniqueness of the above mentioned the martingale problem stated in Proposition \ref{P47}, $X^i$, $i=1,2$, under $\Q^i$ have the same law. 
 Hence, for every Borel set $B \in \mathcal B(C([0,T]))$, we have 
$$
\P^1(X^1 \in B) = \int_\Omega \frac{1}{V_T^1(X^1)} \one_{X^1 \in B}d \Q^1= \int_\Omega \frac{1}{V_T^2(X^2)} \one_{X^2 \in B}d \Q^2= \P^2(X^2 \in B).
$$
Therefore, $X^1$ under $\P^1$ has  the same law as $X^2$ under $\P^2$. Finally, uniqueness holds for the martingale problem in Definition \ref{D:mtpb_hom_Mark}  related to    $\mathcal D_{\mathcal L}$ in \eqref{barD},   $\mathcal L$ in \eqref{Ldistrib} and  $x_0 \in \R$.
\endproof

\section{Further continuity properties}\label{S:cont prop}
We introduce here some continuity properties which are used in the companion paper \cite{BandiniRusso_RevisedWeakDir}.

Let $C_{BUC}(D_{-}(0,\,T); B(0,T))$ be the set of functions  $G: D_{-}(0,T)\rightarrow B(0,T)$ bounded and uniformly continuous on closed balls $B_M \subset D_{-}(0,T)$ of radius $M$.
$C_{BUC}(D_{-}(0,\,T); B(0,T))$ is a Fr\'echet space equipped  with the distance generated by the seminorms
$$
\sup_{ \eta \in B_M} ||G(\eta)||_\infty, \quad M \in \N.
$$
For $f \in C^{1+ \alpha}_{\textup{loc}} \cap  C^0_b$, we  
set  
\begin{equation}\label{E:F}
F^f(y):= \int_{\R} (f(y + x) -f(y)
-k(x)\,f'(y)) Q(y,dx), \quad y \in \R. 
\end{equation}
\begin{proposition}\label{L:Hverified}
Let $\alpha \in [0,1]$.
Assume Hypotheses  \ref{H:h}, \ref{H:hBIS} and \ref{H:h2} with respect to $\alpha$.
Let $H$ be a functional satisfying Hypothesis \ref{H:Hcont} 
and $Q(\cdot, dx)$ be a transition kernel satisfying Hypothesis 
\ref{H:totalvar}    with respect to  $\alpha$.  Assume moreover that $H$ is uniformly continuous on closed balls. Below we will make use of $\mathcal{D}_{\mathcal L}$ and  $\mathcal L$ defined in \eqref{barD} and \eqref{Ldistrib}, respectively.
Then the following holds.
\begin{enumerate}
	\item For every $f \in \mathcal{D}_{\mathcal L}$,  $\mathcal L f\in  C_{BUC}(D_{-}(0,T); B(0,T))$.
	\item The linear map $\mathcal L:  \mathcal{D}_{\mathcal L}\rightarrow C_{BUC}(D_{-}(0,T); B(0,T))$  is continuous.
\end{enumerate}
\end{proposition}
\proof
Let us start by proving item 1.
Let 
$f \in \mathcal D_{\mathcal L}$. 
Let us first show that 
$
\eta\mapsto J^f(\eta)(t):= F^f(\eta(t))
$
belongs to $C_{BUC}(D_{-}(0,T); B(0,T))$.
Let $M>0$. We show  that $J^f$ is bounded and uniformly continuous on $B_M:=\{\eta \in D_{-}(0,T):\,\,||\eta||_\infty \leq M\}$.

Since $F^f$ is continuous by Lemma \ref{L:intermediate}, it is a bounded function on bounded intervals. 
Therefore $J^f$ is bounded, being $B_M$  bounded. 

Let $\delta>0$ and  $\eta_1, \eta_2 \in B_M$ such that 
$$
\sup_{t \in [0,T]}|\eta_1(t)-\eta_2(t)|< \delta.
$$
Then, for every $t \in [0,T]$, 
$$
|F^f(\eta_1(t))-F^f(\eta_2(t))|\leq \sup_{\underset{|y_1| \leq M, |y_2| \leq M, |y_1-y_2|< \delta}{y_1, y_2}}|F^f(y_1)-F^f(y_2)|.
$$
This  implies that $J^f$ is uniformly continuous on $B_M$, since $F^f$ is unifomly continuous on compact sets.

 The map
$
\eta \mapsto  I^f(\eta)(t) :=L f(\eta(t))
$
is bounded and uniformly continuous on $B_M$ because 
 $y \mapsto L f(y)$ is bounded and uniformly continuous on compact intervals.
 
It remains to prove that 
$
\eta \mapsto  W^f(\eta)(t) := H (\eta)(t) (\sigma f')(\eta(t)),
$
is bounded and uniformly continuous on $B_M$. 
The map $\eta\mapsto (\sigma f')(\eta(\cdot))$ is bounded and uniformly continuous by the same reasons as before, while the map $\eta \mapsto H(\eta)$ is bounded and uniformly continuous by assumption.

 Let us now prove item 2. We recall that, for every $f \in \mathcal D_{\mathcal L}$, 
\begin{equation*}
f \mapsto (\mathcal L f)(\eta)(t) = L f(\eta(t)) + \sigma(\eta(t))H(\eta)(t) f'(\eta(t))
+  F^f(\eta(t)), \quad \eta \in D_{-}(0,T), t \in [0,T].
\end{equation*}
Since $\mathcal L$ is the sum of three linear operators, it will be enough to study the continuity at zero. 
We suppose first that 
$
f \mapsto   F^f
$
is continuous from $\mathcal D_{\mathcal L}$ to $C^0$.
This would imply that
$
f \mapsto J^f(\eta) (t)
$
is continuous.
Indeed, let $M>0$ and  $B_M$ be  the closed ball of $D_{-}(0,T)$ with radius $M$. 
We have 
$$
\sup_{t \in [0,T], \eta \in B_M} |J^{f}(\eta) (t)| \leq  \sup_{y: |y|\leq M} |F^{f}(y)|. 
$$

Let us thus prove that $
f \mapsto F^f
$
is continuous from $\mathcal D_{\mathcal L}$ to $C^0$.
For any $f \in  \mathcal{D}_{\mathcal L}$, we decompose $F^f=F^f_1 + F^f_2$, where $F^f_1, F^f_2$ are the functions introduced in \eqref{decF}, namely, for every $y \in \R$, 
\begin{align*}
	F_1^{f}(y)&= \int_{\mathcal B} (f(y + x) -f(y)
-k(x)\,f'(y)) Q(y,dx),\\
	F_2^{f}(y)&=\int_{\R \setminus \mathcal B} (f(y + x) -f(y)
-k(x)\,f'(y)) Q(y,dx),
\end{align*}
and $\mathcal B=[-R, R]$ is a neighborhood  of $x=0$,  such that $k(x) =x$ on $\mathcal B$. 
   We have 
\begin{align*}
	F_1^{f}(y)&= \int_{\mathcal B} (f(y + x) -f(y)
-x\,f'(y)) Q(y,dx)= \int_{\mathcal B} G^{f}(y,x)\,|x|^{1+\alpha} \,Q(y,dx),
\end{align*}
with 
\begin{align}\label{Gf}
G^f(y,x): =\int_0^1 \frac{ f'(y + ax)
-f'(y)  }{|x|^{\alpha}}\,da.
\end{align}
Using that 
\begin{align}\label{est_Btilde}
\sup_{ y \in K, x \in \mathcal  B} G^f(y,x) \leq ||f'||_{\alpha,M+R},
\end{align}
where $||\cdot||_{\alpha,M+R}$ was defined in \eqref{normalphaR},  we get  
$$
\sup_{y : |y|\leq M}|F_1^f(y)| \leq ||f'||_{\alpha, M+R}\sup_{y : |y|\leq M}
||\one_{\mathcal B}(x)|x|^{1+\alpha} Q(y, dx)||_{var},
$$
where previous supremum is finite by Lemma \ref{L:techResult}-b) with $\ell_1(x)= \one_{\mathcal B}(x)$, taking into account Hypothesis \ref{H:totalvar}.
Therefore
 this converges to zero when $f$ converges to zero in $\mathcal D_{\mathcal L}$. This establishes the continuity of $f \mapsto F_1^f$. 
 
 On the other hand,   the continuity of $f \mapsto F_2^f$ follows from the inequality 
\begin{align*}
\sup_{y : |y|\leq M}|F_2^f(y)|\leq  \Big(2 ||f||_\infty + ||k||_\infty \sup_{y : |y|\leq M}|f'(y)|\Big )\sup_{y : |y|\leq M}||\one_{\mathcal B^c}(x)Q(y, dx)||_{var},
\end{align*}
where previous supremum is finite taking into account again  Lemma \ref{L:techResult}-b) with $\ell_2= \one_{\mathcal B^c}(x)$, again  taking into account Hypothesis \ref{H:totalvar}.

We then remark that
$
f \mapsto  I^f(\eta)(t)
$
is continuous. As a matter of fact,  
$$
\sup_{t \in [0,T], \eta \in B_M} |I^{f}(\eta) (t)| \leq  \sup_{y: |y|\leq M} |L f(y)|, 
$$
and this converges to zero when $f$ converges to zero in $\mathcal D_{\mathcal L}$ (and therefore in $\mathcal D_L $),  taking into account the continuity of $L$ by Definition \ref{exc}.

Finally, the continuity of 
$
f \mapsto  W^f(\eta)(t)
$
follows from the boundedness of $H$, and  the fact that, since $f$ converges to zero on $\mathcal D_{\mathcal L}$, then $f'$ converges to zero uniformly on compacts. 
\endproof

\section{New results on  Dirichlet processes}\label{S:Dirichlet}


For a weak Dirichlet process $X$, we will denote by  $X^c$ its unique martingale component, see Proposition 3.2 in \cite{BandiniRusso_RevisedWeakDir}. We start by stating the following result.
\begin{lemma}\label{L:3.2}
	Let $X$ be a Dirichlet process. Then $X$ is a special weak Dirichlet process, and 
	$$
	[X, X]^c=[X^c, X^c].
	$$
\end{lemma}
\begin{remark}\label{R:weakDirnew}
	A   special semimartingale $Y=M+V$ is a Dirichlet process if and only if $V$ is a continuous process. Indeed, $[V, V]= \sum_{s \leq \cdot} |\Delta V_s|^2$. 
\end{remark}
\proof 
It is a consequence of Proposition 5.7 and  Corollary 5.8-(ii)  in \cite{BandiniRusso1},  where $M^c= X^c$ by the uniqueness of the decomposition in Proposition 3.2 in \cite{BandiniRusso_RevisedWeakDir}.
\endproof
We say that \emph{$\nu^X$ does not jump} if
	\begin{equation}\label{R1bis}
 \nu^X(\{t\}\times B)=0\quad \forall t \in [0,T], \quad B \in \mathcal B(\R^\ast).
	\end{equation}
	\begin{remark}
		If \eqref{R1bis} holds true, then obviously
	\begin{equation}\label{R1}
\int_\R x \,\nu^X(\{t\}\times dx)=0, \quad t \in [0,T]. 
	\end{equation}
	The converse is not true. Indeed, consider for instance the case $\nu^X(dt\,dx)= Q(dx) d \psi_t$, with $\psi$ an increasing c\`adl\`ag discontinuous  function and  $Q(dx)$ a symmetric measure, i.e. such that $Q(B) = Q(-B)$, $B \in \mathcal B(\R)$. 
	\end{remark}
	\begin{proposition}\label{R:R}
			If $X$ is a Dirichlet process, then \eqref{R1} holds true. 
\end{proposition}
\proof 		
Suppose that $X$ is a Dirichlet process. Then by 
Lemma \ref{L:3.2}, 
$X$  is
a special weak Dirichlet process and by Corollary 3.22-(ii) in
\cite{BandiniRusso_RevisedWeakDir}, 
	\begin{equation}\label{R2}
		X= X^c + M^{d,X}+ \Gamma^X,
	\end{equation}
	with $X^c$ the unique continuous  martingale part of $X$,  $M^{d, X}=x\,\star (\mu^X- \nu^X)$ and $\Gamma^Y$  a predictable   and  $\F$-orthogonal process. 
	We have therefore
\begin{align}\label{Delta Gamma}
		\Delta \Gamma^X_t &= \int_{\R} x \,\nu^X(\{t\}\times\,dx), \quad t \in [0,T]. 
	\end{align}
	By uniqueness of decomposition of Dirichlet processes and \eqref{R2}, $[\Gamma^X, \Gamma^X]=0$, therefore $\Delta \Gamma^X_t=0$ for all $t \in [0,T]$, and so \eqref{R1} holds true. 
	\endproof 
Let  $\varphi:\R \rightarrow \R$,    $X$ be a c\`adl\`ag process with jump measure $\mu^X$ such that $Y=\varphi(X)$ is a weak Dirichlet process.
We recall that $Y$ is a special weak Dirichlet process  if  and only if there exists a constant $a >0$  such that 
\begin{equation}\label{intY}
	 (\varphi(X_{s-}+x)-\varphi(X_{s-})) \,\one_{\{|x| >a\}} \star \mu^X \in \mathcal{A}^+_{\textup{loc}}, 
	\end{equation}
see Theorem 3.16 in \cite{BandiniRusso_RevisedWeakDir}.
\begin{remark}\label{R:Rb}
The converse of  Proposition \ref{R:R} is not true in general. Indeed, by Remark \ref{R:3.14}, there exist processes $X$ such that $\nu^X$ does not jump (therefore satisfying \eqref{R1})  that nevertheless are not Dirichlet processes, because \eqref{intY} with $\varphi \equiv Id$ is not verified. 
\end{remark}

Suppose $X$ to be a Dirichlet process and $\varphi \in C^1(\R)$. Is $Y=\varphi(X)$ necessarily a Dirichlet process? 

 When $X$ is a continuous  Dirichlet process and $\varphi \in C^1$, then $Y$ is a Dirichlet process, see the proof of Proposition 4.6 of \cite{rvw}. By Lemma \ref{L:3.2}, $Y$ is also a special weak Dirichlet process. Below we discuss the case when $X$ is a discontinuous Dirichlet process. 



\begin{theorem}\label{P:last}
Let $\varphi: \R \rightarrow \R$ be a $C^1$ function and $X$ be a  Dirichlet process. Then  $\varphi(X)$ is a Dirichlet process if and only if   \eqref{intY} holds true for some $a >0$ and 
\begin{equation}\label{R1varphi}
		\int_\R (\varphi(X_{t-}+x) - \varphi(X_{t-})) \,\nu^X(\{t\}\times dx)=0\quad \forall t \in [0,T].
	\end{equation}
\end{theorem}
\begin{remark}
 If $X$ is continuous, then \eqref{intY}  and 
\eqref{R1varphi} are obviously verified, so we retrieve the result stated in the continuous case. 
\end{remark}

\noindent \emph{Proof of Theorem \ref{P:last}.}
Let us set $Y:=\varphi(X)$.
We first prove the direct implication.  
By Lemma \ref{L:3.2},  $Y$ is a special weak Dirichlet process. On the other hand,  by Theorem 3.16  in \cite{BandiniRusso_RevisedWeakDir}, this implies that \eqref{intY} holds for all $a>0$.  Finally, \eqref{R1varphi} follows  from Proposition \ref{R:R} applied to the process  $Y$, since
 $$
 0=\int_{\R} y \,\nu^Y(\{t\}\times\,dy) = \int_\R (\varphi(X_{t-}+x) - \varphi(X_{t-})) \,\nu^X(\{t\}\times dx), \quad t \in[0,T].
 $$

 We prove now the converse implication.  Since $X$ is  a weak Dirichlet process with finite quadratic variation and taking into account   \eqref{intY},
 we can apply Corollary 3.37 in \cite{BandiniRusso_RevisedWeakDir}. Therefore
 $Y$ is a special weak Dirichlet process with decomposition 
 \begin{equation}\label{chainrule_onlyjumps}
Y=  Y_0 + \int_0^\cdot \varphi'(X_s) d X_s^c+   M^{d, \varphi} +\Gamma(\varphi), 
\end{equation}
with $M^{d, \varphi}=(\varphi(s,X_{s-} + x)-\varphi(s,X_{s-}))\,\star (\mu^X- \nu^X)$ and $\Gamma(\varphi)$   predictable   and  $\F$-orthogonal. To show that $Y$ is a Dirichlet process, we need to prove that 
\begin{equation}\label{R6}
\Gamma(\varphi) := Y - Y^c - M^{d, \varphi}
\end{equation}
is a zero quadratic variation process. By \eqref{R6}, we get 
\begin{align}\label{R7}
	[\Gamma(\varphi),\Gamma(\varphi)]&= [Y,Y]+ [Y^c,Y^c]+ [M^{d, \varphi},M^{d, \varphi}]-2 [Y, Y^c] -2 [Y, M^{d, \varphi}], 
	\end{align}
        provided the latter covariation exists.
        In fact we have used that $[Y^c, M^{d, \varphi}]=0$ being $M^{d, \varphi}$  martingale orthogonal.
	Since $M^{d, \varphi} + \Gamma(\varphi)$ is orthogonal, $[Y, Y^c]=[Y^c, Y^c]$. By Proposition 5.3 in \cite{BandiniRusso1},  we get 
	$$
	[M^{d, \varphi},M^{d, \varphi}]= \sum_{s \leq \cdot}|\Delta M^{d, \varphi}_s|^2.
	$$
Collecting previous considerations, \eqref{R7} reads 
\begin{align}\label{R9}
	[\Gamma(\varphi),\Gamma(\varphi)]&= [Y,Y]-[Y^c,Y^c]+\sum_{s \leq \cdot} |\Delta M^{d, \varphi}_s|^2 -2 [Y, M^{d, \varphi}].
	\end{align}
	Provided the latter covariation exists,  $\Gamma(\varphi)$ is a finite quadratic variation process. 
	
	Now, 
	$
	X= X^c + M^{d, X}+ \Gamma^X,
	$
        and $X$ is a Dirichlet process. Therefore, by uniqueness
        of the decomposition of such a process, 
        $\Gamma^X$  is a zero quadratic variation process
and in particular  continuous. 
Therefore we get
	\begin{align*}
		[X, M^{d, \varphi}]= [X^c + M^{d, X}+ \Gamma^X, M^{d, \varphi}]= [M^{d, X}, M^{d, \varphi}]+[\Gamma^X, M^{d, \varphi}],  
	\end{align*}
	and the latter covariation above vanishes since 
	$$
	|[\Gamma^X, M^{d, \varphi}]|\leq \{[\Gamma^X, \Gamma^X][M^{d, \varphi}, M^{d, \varphi}]\}^{1/2}=0. 
	$$
	 So 
	$$
	[X, M^{d, \varphi}]= \sum_{s \leq \cdot} \Delta M^{d, X}_s \Delta M_s^{d, \varphi} =\sum_{s \leq \cdot} \Delta X_s \Delta M_s^{d, \varphi}
	$$
	by Proposition 5.3 in \cite{BandiniRusso1}. It follows that $(X, M^{d, \varphi})$ has all its mutual covariations. 
	Therefore, by Lemma \ref{L:app1}-2, 
	\begin{equation}\label{toen}
	[Y, M^{d, \varphi}] = \int_0^\cdot \varphi'(X_{s-}) d [X, M^{d, \varphi}]_s= \sum_{s \leq \cdot}\varphi'(X_{s-}) \Delta X_s \Delta  \varphi(X_s)
	\end{equation}
	so $[Y, M^{d, \varphi}]$ exists  and, going back to \eqref{R9}, we conclude that $\Gamma(\varphi)$ is a finite quadratic variation process.

	In particular, formula \eqref{toen} gives
	\begin{equation}\label{R10}
		[Y, M^{d, \varphi}]^c =0.
	\end{equation}
	Taking the continuous component in the equality \eqref{R9} and formula \eqref{R10}, we get 
	$$
	[\Gamma(\varphi), \Gamma(\varphi)]^c = [Y, Y]^c- [Y^c, Y^c].
	$$
	By  Lemma \ref{L:app1}-1 and Lemma \ref{L:3.2}, we get 
	$$
	[Y, Y]^c = \int_0^\cdot |\varphi'(X_{s-})|^2 d[X, X]^c_s= \int_0^\cdot |\varphi'(X_{s-})|^2 d[X^c, X^c]_s.
	$$
	By Theorem 3.36 in \cite{BandiniRusso_RevisedWeakDir}, 
	$$
	[Y^c, Y^c]= \int_0^\cdot |\varphi'(X_{s-})|^2 d[X^c, X^c]_s,
	$$
	which implies $[\Gamma(\varphi), \Gamma(\varphi)]^c=0$. It remains to prove that $\Delta \Gamma (\varphi)=0$, since $[\Gamma(\varphi), \Gamma(\varphi)]= [\Gamma(\varphi), \Gamma(\varphi)]^c + \sum_{s \leq \cdot} |\Delta \Gamma(\varphi)|^2$.
By \eqref{R6},
	\begin{equation*}\label{RA0}
		\Delta \Gamma_s(\varphi) = \int_\R (\varphi(X_{s-} +x)- \varphi(X_{s-}))  \nu^X(\{s\}\times dx),
	\end{equation*}
        which is zero by  assumption
  \eqref{R1varphi}.  
\qed

\begin{remark} \label{R316}
\begin{itemize}
\item [a)]	If $\nu^X$ does not jump then obviously
  \eqref{R1varphi} holds true. In this case, according to Theorem \ref{P:last}, $Y=\varphi(X)$ is a Dirichlet process if and only if   \eqref{intY} holds true for some $a >0$.

\item [b)] It is  possible to have a Dirichlet process $X$ and a process $Y= \varphi(X)$, with $\varphi \in C^1$,  that is not a  Dirichlet process. We can indeed show the existence of a martingale $X$ such that $\varphi(X)$ is not even a special weak Dirichlet process:  we will show that \eqref{intY} is not verified, and so, by the direct implication of Theorem  \ref{P:last},
 $Y$ cannot be a Dirichlet process. 
	
To this end, let $Z$ be a Cauchy random variable, in particular its density is
$$
p(x) = \frac{1}{\pi(1+x^2)}, \quad x \in \R. 
$$
We set 
$$
\tilde Z=\sqrt{Z}\one_{\{Z >0\}}+\sqrt{-Z}\one_{\{Z <0\}}.
$$
Clearly, $\E[|\tilde Z|] < \infty$ and $\E[\tilde Z]=0$. We define now 
\begin{align*}
	X_t:=\left\{
	\begin{array}{ll}
	0 &\textup{if}\,\, t \in [0,1[\\
	\tilde Z  &\textup{if}\,\, t >1.
	\end{array}
	\right.
	\end{align*}
	We consider the filtration $\F=(\mathcal F_t)$, with $\mathcal F_t$ being the trivial $\sigma$-algebra for $t \in [0,1[$ and being $\sigma(Z)$ for $t \geq 1$. It follows that $X$ is a martingale: in fact $X_t \in L^1$ for all $t \geq 0$, and 
	\begin{align*}
	\E[X_t|\mathcal F_s]=\left\{
	\begin{array}{ll}
	0=X_s  & \textup{if}\,\, s,t \in [0,1[\\
	\E[X_t]= \E[\tilde Z]=0=X_s  &\textup{if}\,\, s <1, t >1\\
	\tilde Z = X_s & \textup{if}\,\, t >s>1.
	\end{array}
	\right.
	\end{align*}
	 On the other hand, setting $\varphi(x) = x^2$, we have 
	\begin{align*}
	\varphi(X)=\left\{
	\begin{array}{ll}
	0 &\textup{if}\,\, t \in [0,1[\\
	|Z| &\textup{if}\,\, t \geq 1.
	\end{array}
	\right.
	\end{align*}
	Now $Y= \varphi(X)$ is a semimartingale since it is an increasing process, but \eqref{intY} is not verified. As a matter of fact, 
	\begin{align*}
	 (\varphi(X_{s-}+x)-\varphi(X_{s-})) \,\one_{\{|x| >1\}} \star \mu^X_t &= \sum_{s \leq t} \Delta \varphi(X_s)\one_{\{|\Delta X_s| >1\}}=  |Z|\,\one_{\{|Z| >1\}} \one_{\{t \geq 1\}}\notin \mathcal{A}^+_{\textup{loc}},
	\end{align*}
	since $\E[|Z|] = \infty$.
\end{itemize}
\end{remark}

\begin{appendix}
\section{Some technical results}

\label{ATechnical}

\begin{lemma}\label{L:techResult}
 Let  $\ell_1, \ell_2: \R \rightarrow \R$ be Borel measurable and bounded functions such that  $\ell_1$ has  compact support and  $\ell_2$ has support in $\R^\ast$.  
 Set
\begin{align}
	y &\mapsto \ell_1(x) |x|^{1+\alpha} \, Q(y,dx)=:\tilde Q^{\ell_1}(y,dx),\label{Qtildel1}\\
	y &\mapsto  \ell_2(x)  \, Q(y,dx)=:\tilde Q^{\ell_2}(y,dx).\label{Qtildel2}
\end{align}
\begin{itemize}
	\item [a)]If $Q$ satisfies  Hypothesis \ref{H:Kmeas}  for some $\alpha \in [0,\,1]$, then \eqref{Qtildel1}-\eqref{Qtildel2} are   bounded in the total variation norm.
	\item [b)]If $Q$ satisfies Hypothesis \ref{H:totalvar} for some $\alpha \in [0,\,1]$, then  \eqref{Qtildel1}-\eqref{Qtildel2} are  continuous  in the total variation norm.
\end{itemize}
\end{lemma}
\proof
Let $R_1>1$ such that  $\mathcal B_{R_1}:=\{x \in \R: |x|\leq R_1\}$ contains the compact support of $\ell_1$, and  $0 <R_2 \leq 1$ such that   $\mathcal B_{R_2}^c:=\{x \in \R: |x|> R_2\}$ contains the support of  $\ell_2$. Let $M_{\ell_1}:= \sup |\ell_1|$, $M_{\ell_2}:= \sup |\ell_2|$.

a) We first prove that  $\tilde Q^{\ell_1}(\cdot,dx)$ and  $\tilde Q^{\ell_2}(\cdot,dx)$ in \eqref{Qtildel1}-\eqref{Qtildel2} are bounded in the total variation norm. For $y \in \R$, we have 
\begin{align*}
	 \tilde Q^{\ell_1}(y,dx)&= \one_{\{|x| \leq R_1\}}\,\ell_1(x) |x|^{1+\alpha} \, Q(y,dx)\\
	 &=
	 \one_{ \{|x| \leq 1\}}\,\ell_1(x) |x|^{1+\alpha} \, Q(y,dx)+ \one_{\{1 <|x| \leq R_1\}}
	 \,\ell_1(x) |x|^{1+\alpha} \, Q(y,dx).
\end{align*}
We get
\begin{align*}
	 ||\tilde Q^{\ell_1}(y,dx)||_{var}&= ||\one_{\{|x| \leq 1\}}\,\ell_1(x) (|x|^{1+\alpha}\wedge 1) \, Q(y,dx)+ \one_{\{1 <|x| \leq R_1\}}\,\ell_1(x) |x|^{1+\alpha} \, Q(y,dx)||_{var}\\
	 &\leq M_{\ell_1}[||(|x|^{1+\alpha}\wedge 1) \, Q(y,dx)||_{var}+ R_1^{1+\alpha}||\one_{\{1 <|x| \leq R_1\}} \, Q(y,dx)||_{var}]\\
	 &\leq M_{\ell_1}(1+ R_1^{1+\alpha})
	||(1 \wedge |x|^{1+\alpha}) \, Q(y,dx)||_{var}.
\end{align*}
On the other hand, 
\begin{align*}
	\tilde Q^{\ell_2}(y,dx)&=\one_{\{|x|>R_2\}}\ell_2(x)  \, Q(y,dx)=
	\one_{\{R_2 <|x| \leq 1\}}\ell_2(x)  \, Q(y,dx)+
	\one_{\{|x| > 1\}}\ell_2(x)  \, Q(y,dx),
\end{align*}
so that
\begin{align*}
	 ||\tilde Q^{\ell_2}(y,dx)||_{var}&= ||\one_{\{R_2 <|x| \leq 1\}}\,\ell_2(x)  \, Q(y,dx)+ \one_{\{|x|>1}\,\ell_2(x) (|x|^{1+\alpha\}} \wedge 1) \, Q(y,dx)||_{var}\\
	 & \leq \frac{1}{R_2^{1+\alpha}}||\one_{\{R_2 <|x| \leq 1\}}\,\ell_2(x)  \, (|x|^{1+\alpha} \wedge 1)Q(y,dx)||_{var}+ ||\one_{\{|x|>1\}}\,\ell_2(x) (|x|^{1+\alpha} \wedge 1) \, Q(y,dx)||_{var}\\
	  &\leq M_{\ell_2}\Big(1+ \frac{1}{R_2^{1+\alpha}}\Big)
	 ||(1 \wedge |x|^{1+\alpha}) \, Q(y,dx)||_{var}.
 \end{align*}
By Hypothesis \ref{H:Kmeas} together with Remark \ref{R:boundednessTotalVar}, this proves that  $\tilde Q^{\ell_1}(\cdot,dx)$ and  $\tilde Q^{\ell_2}(\cdot,dx)$ in \eqref{Qtildel1}-\eqref{Qtildel2} are bounded in the total variation norm.

b) Let us now prove that $y \mapsto\tilde Q^{\ell_1}(y,dx)$ and  $y \mapsto\tilde Q^{\ell_2}(y,dx)$ in \eqref{Qtildel1}-\eqref{Qtildel2} are continuous in the total variation norm. Let $(y_n)$ be a real sequence converging to $y_0$. 
 We have 
\begin{align*}
	 \tilde Q^{\ell_1}(y_n,dx)-\tilde Q^{\ell_1}(y_0,dx)
	 &= 
	 \one_{\{ |x| \leq 1\}}\,\ell_1(x) |x|^{1+\alpha} \, [Q(y_n,dx)-Q(y_0,dx)]\\
	 &+ \one_{\{1 <|x| \leq R_1\}}
	 \,\ell_1(x) |x|^{1+\alpha} \,[Q(y_n,dx)-Q(y_0,dx)],
\end{align*}
and thus
\begin{align*}
	 ||\tilde Q^{\ell_1}(y_n,dx)-\tilde Q^{\ell_1}(y_0,dx)||_{var}&\leq ||\one_{\{|x| \leq 1\}}\,\ell_1(x) (1 \wedge |x|^{1+\alpha}) \, [Q(y_n,dx)-Q(y_0,dx)||_{var}\\
	 &+ ||\one_{\{1 <|x| \leq R_1\}}\,\ell_1(x) |x|^{1+\alpha} \, [Q(y_n,dx)-Q(y_0,dx)]||_{var}\\
	 &\leq M_{\ell_1}(1+ R_1^{1+\alpha})
	 ||(1 \wedge |x|^{1+\alpha}) (Q(y_n,dx)- Q(y_0,dx))||_{var}.
\end{align*}
On the other hand,  
\begin{align*}
	\tilde Q^{\ell_2}(y_n,dx)-\tilde Q^{\ell_2}(y_0,dx)
	&=
	\one_{\{R_2 < |x| \leq 1\}}\ell_2(x)  \, [Q(y_n,dx)-Q(y_0,dx)]\\
	&+
	\one_{\{|x| > 1\}}\ell_2(x)  \, [Q(y_n,dx)-Q(y_0,dx)],
\end{align*}
so that
\begin{align*}
	 ||\tilde Q^{\ell_2}(y_n,dx)-\tilde Q^{\ell_2}(y_0,dx)||_{var}
	 &\leq ||\one_{\{R_2 <|x| \leq 1\}}\,\ell_2(x)  \, [Q(y_n,dx)-Q(y_0,dx)]||_{var}\\
	 &+ ||\one_{\{|x|>1\}}\,\ell_2(x) (1 \wedge |x|^{1+\alpha}) \, [Q(y_n,dx)-Q(y_0,dx)]||_{var}\\
	  &\leq M_{\ell_2}\big(\frac{1}{R_2^{1+\alpha}}+1\big)
	  ||(1 \wedge |x|^{1+\alpha}) (Q(y_n,dx)-Q(y_0,dx))||_{var}.
\end{align*}
By Hypothesis \ref{H:totalvar}, this proves that  $y \mapsto\tilde Q^{\ell_1}(y,dx)$ and  $y \mapsto\tilde Q^{\ell_2}(y,dx)$ in \eqref{Qtildel1}-\eqref{Qtildel2} are continuous in  total variation topology.
\endproof
\begin{lemma}\label{L:intermediate}
Let $Q(\cdot, dx)$ be a transition kernel satisfying  Hypothesis 
\ref{H:totalvar}  for some $\alpha \in [0,\,1]$. 
Then the function $f \mapsto F^f$, $C^{1+ \alpha}_{\textup{loc}} \cap  C^0_b \rightarrow \R$, defined  in \eqref{E:F}
is continuous.
\end{lemma}
\proof
Let $\mathcal B=[-R,R]$ be a neighborhood  of $x=0$,  such that $k(x) =x$ on $\mathcal B$. 
Define 
\begin{align}\label{decF}
F^f(y)&= \int_{\mathcal B} (f(y + x) -f(y)
-k(x)\,f'(y)) Q(y,dx)
\notag\\
&+\int_{\R\setminus \mathcal B} (f(y + x) -f(y)
-k(x)\,f'(y)) Q(y,dx)=: F_1^f(y) + F_2^f(y).
\end{align}
Let   $(y_n)$ be a sequence converging to $y_0$ in $\R$, and let $K=[-M,M]$ be a compact set containing $(y_n)$.
We start by noticing that    
\begin{align*}
 F^f_1(y) = \int_{\mathcal B}	(f(y + x) -f(y)
-x\,f'(y)) Q(y, dx)  = \int_{\mathcal B } G^f(y,x)\,|x|^{1+\alpha}Q(y, dx), 
\end{align*}
with $G^f$ in \eqref{Gf}.
Then 
\begin{align}\label{Fdif}
	F^f_1(y_n) - F^f_1(y_0) &= \int_{\mathcal B} G^f(y_n,x) \,|x|^{\alpha+1}\, Q(y_n,dx)-\int_{\mathcal B} G^f(y_0,x)\,|x|^{\alpha+1}\, Q(y_0,dx)\notag\\
					&=\int_{\mathcal B} [G^f(y_n,x)-G^f(y_0,x)] \,|x|^{\alpha+1}\,Q(y_0,dx)+ \int_{\mathcal B} G^f(y_n,x)\,|x|^{\alpha+1}\,[Q(y_0,dx)-Q(y_n,dx)].
					\end{align}				
	We recall \eqref{est_Btilde}, namely
\begin{align*}
\sup_{ y \in K, x \in \mathcal  B} G^f(y,x) \leq ||f'||_{\alpha,M+R}.
\end{align*}
For every $x \in \R$,  $y \mapsto G^f(y,x)$ is continuous by Lebesgue dominated convergence theorem. 				
Then the first term in the right-hand side of \eqref{Fdif} converges again by the  Lebesgue dominated convergence  theorem, taking into account
that the measure $\tilde Q^{\ell_1}(y,dx):= \one_{\mathcal B}(x)|x|^{1+\alpha} \, Q(y,dx)$ is finite thanks to  Lemma \ref{L:techResult} with $\ell_1(x) = \one_{\mathcal B}(x)$. 
 The convergence of the second term in the right-hand side of \eqref{Fdif} follows by 
the continuity of $y \mapsto \tilde Q^{\ell_1}(y,dx)$ in the total variation topology due to  Lemma \ref{L:techResult}, noting that 
\begin{align}\label{estG}
\Big|\int_{\mathcal B} G^f(y_n,x)\,|x|^{\alpha+1}\,[Q(y_0,dx)-Q(y_n,dx)]\Big|\leq ||f'||_{\alpha, M+R} \, ||\tilde  Q^{\ell_1}(y_0, dx)-\tilde Q^{\ell_1}(y_n, dx)||_{var}.
\end{align}
On the other hand, 
\begin{align}\label{Fdif2}
	F_2^f(y_n) - F_2^f(y_0) 
	&=\int_{\R\setminus \mathcal B} ((f(y_n + x)-f(y_n)) -(f(y_0+x)-f(y_0)) Q(y_0,dx)\notag\\
	&+(f'(y_n) -f'(y_0))\int_{\R\setminus \mathcal B} k(x)  Q(y_0,dx)\notag\\
&+\int_{\R\setminus \mathcal B} (f(y_n + x) -f(y_n) - k(x) f'(y_n) )\,[Q(y_0,dx)-Q(y_n,dx)]\notag\\
&=:I'(y_n, y_0)+I''(y_n, y_0)+I'''(y_n, y_0).
\end{align}
Since $(y_n)$ lives in the compact $K$ 
 and $f$ is bounded, we have 
\begin{align*}
	I'''(y_n, y_0)&
	=\int_{\R} (f(y_n + x) -f(y_n)- f'(y_n) k(x))\one_{{\mathcal B}^c}(x)\,[Q(y_0,dx)-Q(y_n,dx)]\\
	&\leq C ||\tilde Q^{\ell_2}(y_0,dx)-\tilde Q^{\ell_2}(y_n,dx)||_{var}
\end{align*}
with $ \tilde Q^{\ell_2}(y,dx):=\one_{{\mathcal B}^c}(x)  \, Q(y,dx)$, and $C= 2||f||_\infty + ||k||_\infty \sup_{y \in K} |f'(y)| $. The convergence follows from  the continuity of in the total variation topology due to  Lemma \ref{L:techResult}  applied to $\ell_2(x)=\one_{{\mathcal B}^c}(x)$.
	
On the other hand, the convergence of 
$I'(y_n, y_0)$ follows by the Lebesgue dominated convergence theorem, taking into account that the measure $\tilde Q^{\ell_2}(y_0, dx)$ is finite due to Lemma \ref{L:techResult}, and the fact that $f$ bounded. 

Finally, the convergence of 
$I''(y_n, y_0)$ follows because $f'$ is continuous,  taking into account that the measure $k(x)\one_{{\mathcal B}^c}(x) Q(y_0, dx)$ is finite by Lemma \ref{L:techResult} applied to $\ell_2(x) = k(x)\one_{{\mathcal B}^c}(x)$.
\endproof
\section{Stability of finite quadratic variation processes}

\label{AFQV}

The following result was well understood in the context of F\"ollmer's discretizations,  but was never established in the regularization framework. 
\begin{lemma}\label{L:app1}
	\begin{enumerate}
	\item
	Let $Y= \varphi(X)$, where $\varphi: \R \rightarrow \R$ is a $C^1$ function and $X$ is a  c\`adl\`ag process of finite quadratic variation. Then
	\begin{equation*}
		[Y,Y]_t = \int_0^t (\varphi'(X_{s-}^2  d[X, X]_s^c + \sum_{s \leq t}(\Delta \varphi(X_s))^2.
	\end{equation*} In particular,  $Y$ is also a  finite quadratic variation process.
	\item 	Let $Y^1= \varphi(X^1)$ and $Y^2= \phi(X^2)$, where $\varphi$ and $\phi$ are $C^1$ functions and $X^1, X^2$ are   c\`adl\`ag processes such that $(X^1, X^2)$ has all its mutual covariations. Then 
	\begin{equation*}
		[Y^1,Y^2]_t = \int_0^t \varphi'(X^1_s) \phi'(X^2_{s-}) d[X^1, X^2]_s^c + \sum_{s \leq t}\Delta \varphi(X^1_s)\,\Delta \phi(X^2_s).
	\end{equation*}
	\end{enumerate}

\end{lemma}
\proof
1. Let $t \in [0,\,T]$, $\varepsilon \in [0,\,1]$. We expand, for $s \in [0,\,T]$,
\begin{align*}
	\varphi(X_{(s + \varepsilon)\wedge t})-\varphi(X_{s \wedge t})= I_1^\varphi(s, t, \varepsilon) (X_{(s + \varepsilon)\wedge t}-X_{s \wedge t}), 
\end{align*}
where
 \begin{align*}
   I_1^\varphi(s, t, \varepsilon)= \int_0^1 \varphi'(X_{s \wedge t}+
   a (X_{(s + \varepsilon)\wedge t}-X_{s\wedge t}))\,da.
  \end{align*}
Consequently, 
\begin{align*}
	\frac{1}{\varepsilon}(\varphi(X_{(s + \varepsilon)\wedge t})-\varphi(X_{s \wedge t}))^2 &= \frac{1}{\varepsilon}((I_1^\varphi(s, t, \varepsilon))^2-(\varphi'(X_s))^2) (X_{(s + \varepsilon)\wedge t}-X_{s \wedge t})^2 \\
	&+\frac{1}{\varepsilon}(\varphi'(X_s))^2 (X_{(s + \varepsilon)\wedge t}-X_{s \wedge t})^2. 
\end{align*}
Integrating from $0$ to $t$, we get 
\begin{align}\label{J1+J2}
	\frac{1}{\varepsilon}\int_0^t(\varphi(X_{(s + \varepsilon)\wedge t})-\varphi(X_{s}))^2 ds &= \frac{1}{\varepsilon}\int_0^t((I_1^\varphi(s, t, \varepsilon))^2-(\varphi'(X_s)^2) (X_{(s + \varepsilon)\wedge t}-X_{s})^2 ds\notag\\
	&+\frac{1}{\varepsilon}\int_0^t (\varphi'(X_s))^2 (X_{(s + \varepsilon)\wedge t}-X_{s })^2 ds\notag\\
	&= J_1(t, \varepsilon) + J_2(t, \varepsilon).
\end{align}
We notice that, without restriction of generality, passing to a suitable subsequence, we can suppose (with abuse of notation) that
\begin{equation}\label{convbrac}
[X,X]^\varepsilon:=\frac{1}{\varepsilon} \int_0^\cdot \  (X_{(s + \varepsilon)\wedge \cdot}-X_{s})^2\, ds\underset{\varepsilon \rightarrow 0}{\rightarrow}  [X,X], \quad \textup{uniformly a.s.}
\end{equation}
Since $X$ is a finite quadratic variation process, by Lemma A.5 in \cite{BandiniRusso1}, taking into account Definition A.2 and Corollary A.4-2. in \cite{BandiniRusso1}, if $g$ is a c\`adl\`ag process then 
\begin{align}\label{B9}
	\frac{1}{\varepsilon}\int_0^t g_{s} (X_{(s + \varepsilon) \wedge t} -X_{s})^2  ds=\frac{1}{\varepsilon}\int_0^t g_{s-} (X_{(s + \varepsilon) \wedge t} -X_{s})^2  ds \underset{\varepsilon \rightarrow 0}{\rightarrow} \int_0^t g_{s-} d [X,X]_s, \quad \textup{u.c.p.}
\end{align}
Therefore, taking $g_s= (\varphi'(X_s))^2$ in  \eqref{B9}, we get 
\begin{align}\label{J2conv}
J_2(\cdot, \varepsilon) \underset{\varepsilon \rightarrow 0}{\rightarrow} \int_0^\cdot (\varphi'(X_{s-}))^2 d [X,X]_s, \quad  \textup{u.c.p.}	
\end{align}

Next step consists in proving that
\begin{align}\label{toproveJ1}
	J_1(\cdot, \varepsilon) \underset{\varepsilon \rightarrow 0}{\rightarrow} \sum_{s \leq t}\Big[\Big(\int_0^1 \varphi'(X_{s-} + a \Delta X_s) da\Big)^2 - (\varphi'(X_{s-}))^2\Big](\Delta X_s)^2, \quad \textup{u.c.p.}
\end{align}
We fix a realization $\omega \in \Omega$. Proceeding as in the proof of Proposition 2.14 in \cite{BandiniRusso1}, let  $(t_i)$ be an enumeration of all the jumps of $X(\omega)$ in $[0, \, T]$. We have 
$ \sum_i (\Delta X_{t_i}(\omega))^2 < \infty$. 

Let $\gamma >0$ and $N= N(\gamma)$ such that 
\begin{align}\label{estXsquare}
	\sum_{i=N+1}^\infty  (\Delta X_{t_i}(\omega))^2 \leq \gamma^2. 
\end{align}
We introduce 
\begin{align*}
	A(\varepsilon, N)= \bigcup_{i=1}^N ]t_{i}- \varepsilon, t_i],\quad 
	B(\varepsilon, N)= \bigcup_{i=1}^N ]t_{i- 1}, t_{i}- \varepsilon]=[0, \, T]\setminus A(\varepsilon, N). 
\end{align*}
We decompose 
\begin{align}\label{J1}
	J_1(t, \varepsilon) =\frac{1}{\varepsilon} \int_0^t \one_{A(\varepsilon, N)}(s) J_{1 0}(s, t,  \varepsilon) ds +\int_0^t \one_{B(\varepsilon, N)}(s) J_{1 0}(s, t,  \varepsilon) ds=:J_{1A}(t, \varepsilon, N) + J_{1B}(t, \varepsilon, N),
\end{align}
where we have denoted 
$$
J_{1 0}(s, t,  \varepsilon):=(X_{(s + \varepsilon)\wedge t}-X_{s})^2((I_1^\varphi(s, t, \varepsilon))^2-(\varphi'(X_s)^2). 
$$
By Lemma 2.11 in \cite{BandiniRusso1}, it follows that,  uniformly in $t \in [0,\,T]$,
\begin{align}\label{6}
	J_{1A}(t, \varepsilon, N) 
	&\underset{\varepsilon \rightarrow 0}{\rightarrow} \sum_{i=1}^N \one_{]0,\,t]}(t_i)(\Delta X_{t_i})^2\Big(\Big(\int_0^1 \varphi'(X_{t_i -} + a \Delta X_{t_i}) da\Big)^2- (\varphi'(X_{t_i -}))^2\Big).
	\end{align}
On the other hand,
\begin{align*}
	J_{1B}(t, \varepsilon, N) &= \sum_{i=1}^{N}\frac{1}{\varepsilon} \int_0^t \ (X_{(s + \varepsilon)\wedge t}-X_{s})^2 I^{\varphi,i}_{1 B}(s, t,  \varepsilon) \,ds,
\end{align*}  
where 
\begin{align*}
	I^{\varphi,i}_{1B}(s, t,  \varepsilon) &= \one_{]t_{i-1}, t_{i}- \varepsilon]}(s)  \Big[\Big(\int_0^1 \varphi'(X_{s \wedge t}+a(X_{(s + \varepsilon)\wedge t}-X_{s\wedge t}))\,da\Big)^2-(\varphi'(X_s))^2\Big]\\
	& =\one_{]t_{i-1}, t_{i}- \varepsilon]}(s)  \Big[\int_0^1 \varphi'(X_{s \wedge t}+a(X_{(s + \varepsilon)\wedge t}-X_{s\wedge t}))\,da-\varphi'(X_s)\Big]\cdot\\
	&\,\,\,\,\cdot\Big[\int_0^1 \varphi'(X_{s \wedge t}+a(X_{(s + \varepsilon)\wedge t}-X_{s\wedge t}))\,da+\varphi'(X_s)\Big]. 
\end{align*}
For every $i =1,..., N$, we have 
\begin{align*}
	|I^{\varphi,i}_{1 B}(s, t,  \varepsilon)| \leq 2 \sup_{y \in [X_s, X_{s +\varepsilon}]}|\varphi'(y)|\,\delta\Big(\varphi, \, \sup_i \sup_{\underset{|p-q|\leq \varepsilon}{p, q \in [t_{i-1}, t_i]}}|X_p - X_q|\Big). 
\end{align*}
 By Lemma 2.12 in \cite{BandiniRusso1}, there is $\varepsilon_0$ such that, if $\varepsilon < \varepsilon_0$, then 
 \begin{align*}
 	|I^{\varphi,i}_{1 B}(s, t,  \varepsilon)| \leq  2 \sup_{y \in [-||X||_\infty, ||X||_\infty]}|\varphi'(y)|\,\delta(\varphi, \,3 \gamma). 
 	\end{align*}
Consequently, for $\varepsilon < \varepsilon_0$,  
\begin{align}\label{9}
	\sup_{t \in [0,\,T]}|J_{1B}(t, \varepsilon, N)| &\leq    2 \sup_{y \in [-||X||_\infty, ||X||_\infty]}|\varphi'(y)|\,\delta(\varphi, \,3 \gamma)\, \sup_{t \in [0,\,T]}[X,X]^\varepsilon_t,
\end{align}
where the latter supremum is finite by  \eqref{convbrac}.
Going back to \eqref{J1} we get  
\begin{align}\label{J1est}
	&\sup_{t \in [0,\,T]}\Big|J_1(t, \varepsilon)-\sum_{i=1}^\infty \one_{]0,\,t]}(t_i)(\Delta X_{t_i})^2\Big[\Big(\int_0^1 \varphi'(X_{t_i-} + a \Delta X_{t_i}) da\Big)^2 - (\varphi'(X_{t_i-}))^2\Big]\Big| \notag\\
	&\leq \sup_{t \in [0,\,T]}\Big|J_{1A}(t, \varepsilon, N)-\sum_{i=1}^N \one_{]0,\,t]}(t_i)(\Delta X_{t_i})^2\Big[\Big(\int_0^1 \varphi'(X_{t_i-} + a \Delta X_{t_i}) da\Big)^2 - (\varphi'(X_{t_i-}))^2\Big]\Big| \notag\\
	&+ \Big|\sum_{i=N+1}^\infty \one_{]0,\,T]}(t_i)(\Delta X_{t_i})^2\Big(\int_0^1 \varphi'(X_{t_i-} + a \Delta X_{t_i}) da\Big)^2 \Big| + \sup_{t \in [0,\,T]}|J_{1B}(t, \varepsilon, N)|.
\end{align}
Taking the $\limsup_{\varepsilon \rightarrow 0}$ in \eqref{J1est}, collecting \eqref{estXsquare}, \eqref{6} and \eqref{9},  we get 
\begin{align*}	
&\limsup_{\varepsilon \rightarrow 0}\sup_{t \in [0,\,T]}\Big|J_1(t, \varepsilon)-\sum_{i=1}^\infty \one_{]0,\,t]}(t_i)(\Delta X_{t_i})^2\Big[\Big(\int_0^1 \varphi'(X_{t_i-} + a \Delta X_{t_i}) da\Big)^2 - (\varphi'(X_{t_i-}))^2\Big]\Big| \notag\\
	&\leq \sum_{i=N+1}^\infty \one_{]0,\,T]}(t_i)(\Delta X_{t_i})^2 \sup_{y \in [-||X||_\infty, ||X||_\infty]}|\varphi'(y)|
	+ 2 \sup_{y \in [-||X||_\infty, ||X||_\infty]}|\varphi'(y)|\,\delta(\varphi, \,3 \gamma)\sup_{\varepsilon < \varepsilon_0}  \sup_{t \in [0,\,T]}[X,X]_t^\varepsilon
	\\
	&\leq \Big[ \gamma^2 +2 \sup_{\varepsilon < \varepsilon_0}  \sup_{t \in [0,\,T]}[X,X]_t^\varepsilon\, \delta(\varphi, \,3 \gamma)\Big] \sup_{y \in [-||X||_\infty, ||X||_\infty]}|\varphi'(y)|.
\end{align*}
Since $\gamma$ is arbitrary and $\varphi'$ is uniformly continuous on compact intervals, then \eqref{toproveJ1} is proved. 
By  \eqref{J2conv} and \eqref{toproveJ1}, and the fact that $[X, X]= [X, X]^c + \sum_{s \leq t}(\Delta X_s)^2$,  \eqref{J1+J2} yields 
\begin{align*}
	\frac{1}{\varepsilon}\int_0^t(\varphi(X_{(s + \varepsilon)\wedge t})-\varphi(X_{s}))^2 ds&\underset{\varepsilon \rightarrow 0}{\rightarrow}\int_0^\cdot (\varphi'(X_{s-}))^2 d [X,X]^c_s+\sum_{s \leq t}(\varphi'(X_{s-}))^2(\Delta X_s)^2\\
	&+\sum_{s \leq t}\Big[\Big(\int_0^1 \varphi'(X_{s-} + a \Delta X_s) da\Big)^2 - (\varphi'(X_{s-}))^2\Big](\Delta X_s)^2\\
	&=\int_0^\cdot (\varphi'(X_{s-}))^2 d [X,X]^c_s+\sum_{s \leq t}\Big(\int_0^1 \varphi'(X_{s-} + a \Delta X_s) da\Big)^2(\Delta X_s)^2, \quad \textup{u.c.p.}
\end{align*}
The result follows because 
$$
\Delta \varphi(X_s) = \varphi(X_s)-\varphi(X_{s-})= \Delta X_s\int_0^1 \varphi'(X_{s-} + a \Delta X_s) da.
$$

 2. The result follows from point 1. by polarity arguments. 
\endproof

\section{Recalls on semimartingales with jumps}

\label{A:sem}

We recall that a special semimartingale is a semimartingale $X$ which admits a decomposition $X=N+V$, where  $N$ is a local martingale and $V$ is a  finite variation and predictable process, see Definition 4.21, Chapter I,  in \cite{JacodBook}.  Fixing $V_0=0$, such a decomposition is unique, and is called canonical decomposition of $X$, see respectively Proposition 3.16 and Definition 4.22, Chapter I, in \cite{JacodBook}. 

Assume now that $X$ is a semimartingale 
with jump measure $\mu^X$.
	   Given $k \in \mathcal K$, the process
 $
 X^k:= X- \sum_{s \leq \cdot} [\Delta X_s - k(\Delta X_s)]
 $
 is a special semimartingale with unique decomposition
\begin{equation}\label{dec_Xk_semimart}
	X^k=  X^c + M^{k,d} + B^{k, X}, 
\end{equation}
where  $M^{k,d}$ is  a purely discontinuous local martingale,   $X^c$ is the unique continuous martingale part of $X$ (it coincides with the process $X^c$ introduced in  Proposition 3.2 in \cite{BandiniRusso_RevisedWeakDir}), and $B^{k, X}$ is a predictable process of bounded variation. 

According to Definition 2.6, Chapter II in \cite{JacodBook}, the characteristics of $X$ associated with $k \in \mathcal  K$
are  then given by the triplet $(B^k,C,\nu)$ on $(\check \Omega, \check {\mathcal F}, \check {\mathbb F})$ such that
\begin{itemize}
	\item[(i)] $B^k$ is  $\check{\mathbb F}$-predictable, with finite variation on finite intervals, and $B_0^k=0$, namely $B^{k, X}= B^k \circ X$ is the process in \eqref{dec_Xk_semimart};
	\item[(ii)] $C$ is 
	a continuous process of finite variation with $C_0=0$, namely $C^X:= C \circ X=  \langle X^c, X^c\rangle$;
	\item[(iii)] $\nu$ is an $\check{\mathbb F}$-predictable random measure on $\R_+ \times \R$, namely    $\nu^X: = \nu \circ X$ is the compensator of  $\mu^X$.
	\end{itemize}

\begin{theorem}[Theorem  2.42, Chapter II, in \cite{JacodBook}]\label{T: equiv_mtgpb_semimart}
Let $X$ be an adapted  c\`adl\`ag process. 	Let $B^k$ be an  $\check{\mathbb F}$-predictable process, with finite variation on finite intervals, and $B_0^k=0$, 
 $C$ be
	an $\check{\mathbb F}$-adapted continuous process of finite variation with $C_0=0$, 	and  $\nu$ be an $\check{\mathbb F}$-predictable random measure on $\R_+ \times \R$.
There is equivalence between the two following statements.
\begin{itemize} \item [(i)] X is a real  semimartingale with characteristics $(B^k, C, \nu)$. \item [(ii)] For each bounded function $f$ of class $C^2$, the process 
\begin{align*}
&f(X_{\cdot}) - f(X_0) - \frac{1}{2} \int_0^{\cdot}  f''(X_s) \,dC^X_s-  \int_0^{\cdot}  f'(X_{s-}) \,d B_s^{k,X}\\
&-  \int_0^{\cdot}\int_{\R} (f(X_{s-} + x) -f(X_{s-})-k(x)\,f'(X_{s-}))\,\nu^X(ds\,dx)
\end{align*} 
is a local martingale.
\end{itemize} 
\end{theorem}
%
%
%
%
%
Let    $(\Omega, \mathcal F, \mathbb F)$ be the canonical filtered space, and   $X$ the canonical process.  
Let moreover $\mathcal H$ be another $\sigma$-algebra and $\P_{\mathcal H}$ be a probability measure on $(\Omega, \mathcal H)$. 	
	

	\begin{definition}[Definition  2.4, Chapter III, in \cite{JacodBook}]\label{D:mrtg_pb_Jacod}
	A solution to the martingale problem associated to $(\mathcal H, X)$ and $(\P_{\mathcal H}; B^k, C, \nu)$ is a probability measure $\P$ on $(\Omega, \mathcal F)$ such that  
	\begin{itemize}
		\item [(i)] the restriction $\P|_{\mathcal H}$ of $\P$ to $\mathcal H$ equals $\P_{\mathcal H}$;
		\item [(ii)]$X$ is a semimartingale on the basis $(\Omega, \mathcal F, \mathbb F, \P)$ with characteristics $(B^k,C,\nu)$. 
	\end{itemize}
	We denote by $s(\mathcal H, X|\P_{\mathcal H}; B^k, C, \nu)$ the set of all solutions $\P$. 
\end{definition}
	

\begin{definition}[Definition  2.18, Chapter III, in \cite{JacodBook}]\label{D:2.18}
  Let $b^k: \R_+ \times \R \rightarrow \R$ be Borel,
  $c: \R_+ \times \R \rightarrow \R$ be Borel and nonnegative, and $Q_s(x,dy)$ be a transition kernel from $(\R \times \R_+, \mathcal{B} (\R \times \R_+)$ into $(\R, \mathcal{B}(\R))$, with $Q_s(x,\{0\})=0$.
	 $X$ is called a   diffusion with jumps on $(\Omega, \mathcal F, \mathbb F, \P)$ related to $(b^k, c, Q_s(x,dy))$ if it is a semimartingale with  characteristics 
	\begin{align}\label{charact}
		B_t^k = \int_0^t b^k(s, \check X_s) ds,\quad C_t = \int_0^t c(s,\check X_s) ds, \quad 
		\nu(ds\,dx)= Q_s(\check X_{s-},dx) ds.
	\end{align}
\end{definition}

\begin{remark}
Suppose that $X$  is as in  Definition \ref{D:2.18}, then $\P$ is a solution to the martingale problem associated to $(\mathcal H, X)$ and $(\P_{\mathcal H}; B^k, C, \nu)$, with $\mathcal H=\sigma(X_0)$ and $\P_{\mathcal H}$  the law of $X_0$.  
\end{remark}

\begin{hypothesis}\label{H_exis_uni_mtgpb}
  Let $b^k: \R_+ \times \R \rightarrow \R$ be Borel, $c: \R_+ \times \R \rightarrow \R$ be Borel and nonnegative, and $Q_s(x,de)$ be a transition kernel from
  $(\R_+ \times \R, \mathcal{B} (\R_+ \times \R))$   into $(\R,\mathcal{B} (\R))$, with $Q_s(x,\{0\})=0$.
We  assume that
\begin{itemize}
\item [(i)]	$b^k$  is bounded;
\item [(ii)]	 $c$ is bounded, continuous  on $\R_+ \times \R$ and not vanishing at zero;
\item [(iii)] the functions $(t,y) \mapsto \int_A (|x|^2 \wedge 1) Q_t(y,dx)$ are bounded and continuous for all $A \in \mathcal B(\R)$. 
\end{itemize}
\end{hypothesis}

\begin{theorem}[Theorem  2.34, Chapter III, in \cite{JacodBook}]	\label{T:2.34}
Let  $(B^k, C, \nu)$ be of the type in \eqref{charact}, and  such that   $b^k(s,x)$, $c(s,x)$, $Q_s(x,dy)$ satisfy Hypothesis \ref{H_exis_uni_mtgpb}. 
Then, there is a transition kernel $P_{z}(d \omega)$ from
$(\R, \mathcal{B} (\R))$ into $(\Omega, \mathcal F)$ with the following property: for every $z$, $P_{z}$ is the unique probability measure under which the canonical process $X$ is a diffusion with jumps, with $P_{z}(X_0 =z)=1$ and with characteristics given by \eqref{charact}.
\end{theorem}

\begin{corollary}
	Let $P_z(d \omega)$ be the transition kernel introduced in Theorem \ref{T:2.34}. Let $z \in \R$. Then $P_z$ is the unique solution to the martingale problem 
associated to $(\mathcal H, X)$ and $(\P_{\mathcal H}; B^k, C, \nu)$, with $\mathcal H = \sigma(X_0)$ and $\P_{\mathcal H}$  determined by $\P(X_0 \in B) = \delta_z(B)$.
\end{corollary}

\end{appendix}

\bibliographystyle{imsart-number} 
\bibliography{../../BIBLIO_FILE/BiblioLivreFRPV_TESI} 

\def\cprime{$'$} \def\cprime{$'$} \def\cprime{$'$} \def\cprime{$'$}
\begin{thebibliography}{20}

\bibitem{mytnik}
\begin{barticle}[author]
\bauthor{\bsnm{Athreya},~\bfnm{S.}\binits{S.}},
  \bauthor{\bsnm{Butkovsky},~\bfnm{O.}\binits{O.}} \AND
  \bauthor{\bsnm{Mytnik},~\bfnm{L.}\binits{L.}}
(\byear{2020}).
\btitle{Strong existence and uniqueness for stable stochastic differential
  equations with distributional drift}.
\bjournal{Ann. Probab.}
\bvolume{48}
\bpages{178--210}.
\bdoi{10.1214/19-AOP1358}
\end{barticle}
\endbibitem

\bibitem{bahouri}
\begin{bbook}[author]
\bauthor{\bsnm{Bahouri},~\bfnm{H.}\binits{H.}},
  \bauthor{\bsnm{Chemin},~\bfnm{J-Y.}\binits{J.-Y.}} \AND
  \bauthor{\bsnm{Danchin},~\bfnm{R.}\binits{R.}}
(\byear{2011}).
\btitle{Fourier Analysis and Nonlinear Partial Differential Equations}.
\bpublisher{Springer}.
\end{bbook}
\endbibitem

\bibitem{BandiniRusso1}
\begin{barticle}[author]
\bauthor{\bsnm{Bandini},~\bfnm{E.}\binits{E.}} \AND
  \bauthor{\bsnm{Russo},~\bfnm{F.}\binits{F.}}
(\byear{2017}).
\btitle{Weak {D}irichlet processes with jumps.}
\bjournal{{Stochastic Process. Appl.}}
\bvolume{127}
\bpages{4139-4189}.
\end{barticle}
\endbibitem

\bibitem{BandiniRusso_RevisedWeakDir}
\begin{barticle}[author]
\bauthor{\bsnm{Bandini},~\bfnm{E.}\binits{E.}} \AND
  \bauthor{\bsnm{Russo},~\bfnm{F.}\binits{F.}}
(\byear{2022}).
\btitle{Weak {D}irichlet processes and generalized martingale problems}.
\bjournal{Preprint Arxiv 2205.03099}.
\end{barticle}
\endbibitem

\bibitem{bq}
\begin{barticle}[author]
\bauthor{\bsnm{Bass},~\bfnm{R.~F.}\binits{R.~F.}} \AND
  \bauthor{\bsnm{Chen},~\bfnm{Z.~Q.}\binits{Z.~Q.}}
(\byear{2001}).
\btitle{Stochastic differential equations for {D}irichlet processes}.
\bjournal{Probab. Theory Related Fields}
\bvolume{121}
\bpages{422--446}.
\end{barticle}
\endbibitem

\bibitem{cannizzaro}
\begin{barticle}[author]
\bauthor{\bsnm{{Cannizzaro}},~\bfnm{G.}\binits{G.}} \AND
  \bauthor{\bsnm{{Chouk}},~\bfnm{K.}\binits{K.}}
(\byear{2018}).
\btitle{{Multidimensional SDEs with singular drift and universal construction
  of the polymer measure with white noise potential}}.
\bjournal{{Ann. Probab.}}
\bvolume{46}
\bpages{1710--1763}.
\end{barticle}
\endbibitem

\bibitem{deraynal2020}
\begin{barticle}[author]
\bauthor{\bparticle{Chaudru~de} \bsnm{Raynal},~\bfnm{P.~E.}\binits{P.~E.}} \AND
  \bauthor{\bsnm{Menozzi},~\bfnm{S.}\binits{S.}}
(\byear{2020}).
\btitle{On multidimensional stable-driven stochastic Differential Equations
  with {B}esov drift}.
\bjournal{Preprint Arxiv 1907.12263}.
\end{barticle}
\endbibitem

\bibitem{diel}
\begin{barticle}[author]
\bauthor{\bsnm{Delarue},~\bfnm{F.}\binits{F.}} \AND
  \bauthor{\bsnm{Diel},~\bfnm{R.}\binits{R.}}
(\byear{2016}).
\btitle{Rough paths and 1d {SDE} with a time dependent distributional drift:
  application to polymers}.
\bjournal{Probab. Theory Related Fields}
\bvolume{165}
\bpages{1--63}.
\end{barticle}
\endbibitem

\bibitem{issoglio}
\begin{barticle}[author]
\bauthor{\bsnm{Flandoli},~\bfnm{F.}\binits{F.}},
  \bauthor{\bsnm{Issoglio},~\bfnm{E.}\binits{E.}} \AND
  \bauthor{\bsnm{Russo},~\bfnm{F.}\binits{F.}}
(\byear{2017}).
\btitle{Multidimensional stochastic differential equations with distributional
  drift}.
\bjournal{Trans. Amer. Math. Soc.}
\bvolume{369}
\bpages{1665--1688}.
\bdoi{10.1090/tran/6729}
\bmrnumber{3581216}
\end{barticle}
\endbibitem

\bibitem{frw1}
\begin{barticle}[author]
\bauthor{\bsnm{Flandoli},~\bfnm{F.}\binits{F.}},
  \bauthor{\bsnm{Russo},~\bfnm{F.}\binits{F.}} \AND
  \bauthor{\bsnm{Wolf},~\bfnm{J.}\binits{J.}}
(\byear{2003}).
\btitle{Some {SDE}s with distributional drift. {I}. {G}eneral calculus}.
\bjournal{Osaka J. Math.}
\bvolume{40}
\bpages{493--542}.
\end{barticle}
\endbibitem

\bibitem{frw2}
\begin{barticle}[author]
\bauthor{\bsnm{Flandoli},~\bfnm{F.}\binits{F.}},
  \bauthor{\bsnm{Russo},~\bfnm{F.}\binits{F.}} \AND
  \bauthor{\bsnm{Wolf},~\bfnm{J.}\binits{J.}}
(\byear{2004}).
\btitle{Some {SDE}s with distributional drift. {II}. {L}yons-{Z}heng structure,
  {I}t\^o's formula and semimartingale characterization}.
\bjournal{Random Oper. Stochastic Equations}
\bvolume{12}
\bpages{145--184}.
\end{barticle}
\endbibitem

\bibitem{chineseBook}
\begin{bbook}[author]
\bauthor{\bsnm{He},~\bfnm{S.}\binits{S.}},
  \bauthor{\bsnm{Wang},~\bfnm{J.}\binits{J.}} \AND
  \bauthor{\bsnm{Yan},~\bfnm{J.}\binits{J.}}
(\byear{1992}).
\btitle{Semimartingale theory and stochastic calculus}.
\bpublisher{Science Press Bejiing New York}.
\end{bbook}
\endbibitem

\bibitem{jacod_book}
\begin{bbook}[author]
\bauthor{\bsnm{Jacod},~\bfnm{J.}\binits{J.}}
(\byear{1979}).
\btitle{Calcul {S}tochastique et {P}robl\`emes de martingales}.
\bseries{Lecture Notes in Mathematics}
\bvolume{714}.
\bpublisher{Springer}, \baddress{Berlin}.
\end{bbook}
\endbibitem

\bibitem{JacodBook}
\begin{bbook}[author]
\bauthor{\bsnm{Jacod},~\bfnm{J.}\binits{J.}} \AND
  \bauthor{\bsnm{Shiryaev},~\bfnm{A.~N.}\binits{A.~N.}}
(\byear{2003}).
\btitle{Limit theorems for stochastic processes},
\bedition{second} ed.
\bseries{Grundlehren der Mathematischen Wissenschaften [Fundamental Principles
  of Mathematical Sciences]}
\bvolume{288}.
\bpublisher{Springer-Verlag, Berlin}.
\end{bbook}
\endbibitem

\bibitem{kremp}
\begin{barticle}[author]
\bauthor{\bsnm{Kremp},~\bfnm{H.}\binits{H.}} \AND
  \bauthor{\bsnm{Perkowski},~\bfnm{N.}\binits{N.}}
(\byear{2022}).
\btitle{Multidimensional {SDE} with distributional drift and {L{\'e}vy} noise}.
\bjournal{Bernoulli}
\bvolume{28}
\bpages{1757--1783}.
\bdoi{10.3150/21-BEJ1394}
\end{barticle}
\endbibitem

\bibitem{ling}
\begin{barticle}[author]
\bauthor{\bsnm{Ling},~\bfnm{C.}\binits{C.}} \AND
  \bauthor{\bsnm{Zhao},~\bfnm{G.}\binits{G.}}
(\byear{2022}).
\btitle{Nonlocal elliptic equation in {H{\"o}lder} space and the martingale
  problem}.
\bjournal{J. Differ. Equations}
\bvolume{314}
\bpages{653--699}.
\bdoi{10.1016/j.jde.2022.01.025}
\end{barticle}
\endbibitem

\bibitem{ORT1_PartI}
\begin{barticle}[author]
\bauthor{\bsnm{Ohashi},~\bfnm{A.}\binits{A.}},
  \bauthor{\bsnm{Russo},~\bfnm{F.}\binits{F.}} \AND
  \bauthor{\bsnm{Teixeira},~\bfnm{A.}\binits{A.}}
(\byear{2022}).
\btitle{On path-dependent {SDEs} involving distributional drifts}.
\bjournal{Mod. Stoch., Theory Appl.}
\bvolume{9}
\bpages{65--87}.
\bdoi{10.15559/21-VMSTA197}
\end{barticle}
\endbibitem

\bibitem{po1}
\begin{bbook}[author]
\bauthor{\bsnm{Portenko},~\bfnm{N.~I.}\binits{N.~I.}}
(\byear{1990}).
\btitle{Generalized diffusion processes}.
\bseries{Translations of Mathematical Monographs}
\bvolume{83}.
\bpublisher{American Mathematical Society}, \baddress{Providence, RI}.
\bnote{Translated from the Russian by H. H. McFaden}.
\end{bbook}
\endbibitem

\bibitem{rtrut}
\begin{barticle}[author]
\bauthor{\bsnm{Russo},~\bfnm{F.}\binits{F.}} \AND
  \bauthor{\bsnm{Trutnau},~\bfnm{G.}\binits{G.}}
(\byear{2007}).
\btitle{Some parabolic {PDE}s whose drift is an irregular random noise in
  space}.
\bjournal{Ann. Probab.}
\bvolume{35}
\bpages{2213--2262}.
\end{barticle}
\endbibitem

\bibitem{rvw}
\begin{barticle}[author]
\bauthor{\bsnm{Russo},~\bfnm{F.}\binits{F.}},
  \bauthor{\bsnm{Vallois},~\bfnm{P.}\binits{P.}} \AND
  \bauthor{\bsnm{Wolf},~\bfnm{J.}\binits{J.}}
(\byear{2001}).
\btitle{A generalized class of {L}yons-{Z}heng processes}.
\bjournal{Bernoulli}
\bvolume{7}
\bpages{363--379}.
\end{barticle}
\endbibitem

\end{thebibliography}


\end{document}